%% file: lamg.tex
\documentclass[final]{siamltex}

\include{styles}
\usepackage{array}

\title{LEAN ALGEBRAIC MULTIGRID (LAMG): \\ FAST GRAPH LAPLACIAN LINEAR SOLVER}
\author{
Oren E. Livne
\thanks{Institute for Genomics \& Systems Biology, The University of Chicago, 900 E 57th St, KCBD 10146, Chicago, IL 60637. Tel: +1-773-702-9765. Fax: +1-773-834-2877. Email: {\tt olivne@uchicago.edu}}
\and
Achi E. Brandt \thanks{Department of Mathematics and Computer Science, The Weizmann Institute of Science, POB 26 Rehovot 76100, Israel. Tel. +972-8-934-3545. Fax: +972-8-934-6023. Email: {\tt abrandt@math.ucla.edu}}
}

\begin{document}
\maketitle

\centerline{\it Dedicated to J. Brahms' Symphony No.~1 in C minor, Op.~68}

\begin{abstract}
Laplacian matrices of graphs arise in large-scale computational applications such as machine learning; spectral clustering of images, genetic data and web pages; transportation network flows; electrical resistor circuits; and elliptic partial differential equations discretized on unstructured grids with finite elements. A Lean Algebraic Multigrid (LAMG) solver of the linear system $Ax=b$ is presented, where $A$ is a graph Laplacian. LAMG's run time and storage are linear in the number of graph edges. LAMG consists of a setup phase, in which a sequence of increasingly-coarser Laplacian systems is constructed, and an iterative solve phase using multigrid cycles. General graphs pose algorithmic challenges not encountered in traditional applications of algebraic multigrid. LAMG combines a lean piecewise-constant interpolation, judicious node aggregation based on a new node proximity definition, and an energy correction of the coarse-level systems. This results in fast convergence and substantial overhead and memory savings. A serial LAMG implementation scaled linearly for a diverse set of \numgraphs real-world graphs with up to \maxedges edges. This multilevel methodology can be fully parallelized and extended to eigenvalue problems and other graph computations.
\end{abstract}



\pagestyle{myheadings}
\thispagestyle{plain}
\markboth{O. E. LIVNE AND A. E. BRANDT}{LEAN ALGEBRAIC MULTIGRID (LAMG) FOR THE GRAPH LAPLACIAN}

\section{Introduction}
\label{introduction}
Let $G=(\cN,\cE,w)$ be a weighted undirected graph, where $\cN$ is a set of $n$ nodes, $\cE$ is a set of $m$ edges, and $w: \cE \rightarrow \Real \backslash \{0\}$ is a weight function. The Laplacian matrix $\bA_{n \times n}$ is naturally defined in terms of its quadratic {\it energy},
\be
    E(\bx) = \bx^{T} \bA \bx = \sum_{(u,v) \in E} w_{uv} \lp x_u-x_v \rp^2\,,\qquad \bx \in \Real^{\cN}\,,
    \label{energy}
\ee
where $T$ denotes the transpose operator. In matrix form,
\be
    \bA = \lp a_{uv} \rp_{u,v}\,,\quad
    a_{uv} :=
    \begin{cases}
        \sum_{v \in \cA_u} w_{uv}\,,& u=v\,, \\
        -w_{uv}\,,&v \in \cA_u\,,
    \end{cases}\quad
    \cA_u := \left\{v: (u,v) \in \cE\right\}\,.
    \label{lap}
\ee
\begin{figure}[h]
\centering
\begin{minipage}{.4\linewidth}
     \includegraphics[height=1.1in]{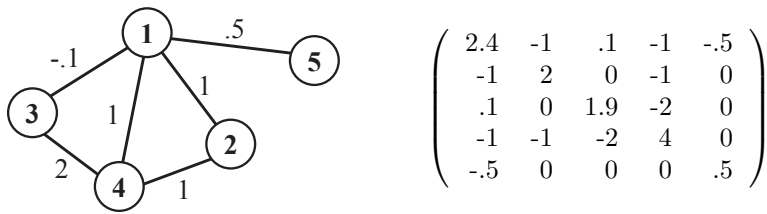}
\end{minipage}
\quad
\begin{minipage}{.4\linewidth}
$\left(
    \begin{tabular}{rrrrr}
   2.4 & -1 & .1 & -1  & -.5 \\
   -1   &  2 &   0 & -1  &    0 \\
    .1 &  0 & 1.9 & -2  &    0 \\
   -1   & -1 &  -2 &  4  &    0 \\
   -.5 &  0 &   0 &  0  &  .5
\end{tabular}
\right)$
\end{minipage}
    \caption{A $5$-node graph with a negative weight. The Laplacian is still semi-positive definite here.}
\label{simple_graph}
\end{figure}

$\bA$ is zero row-sum, Symmetric Positive Semi-definite (SPS), and has $2m+1$ non-zeros. Typically, $m \ll n^2$ and $\bA$ is sparse. Unlike graph-theoretic works \cite{chung, koutis, st06}, $w$ is not assumed to be strictly positive: we also include SPS matrices with positive off-diagonal entries such as high-order and anisotropic grid discretizations \cite{alg_distance_anis}.

$\bA$ is singular: if $G$ has $M$ connected components whose node sets are $\cN_1,...,\cN_M$, then $\bA$ has an eigenvalue $\lambda=0$ of multiplicity $M$ and eigenvectors $\bu_1,\dots,\bu_M$, where $\bu_i$ is the characteristic function of $\cN_i$. In a singly connected graph, the null space is spanned by the vector of ones, $\bu_1=:\bu$. We consider the compatible linear system
\begin{subequations} \label{linsys}
\begin{eqnarray}
    &&\bA \bx = \bb \label{system} \\
    &&\bU^{T} \bx = \balpha\,,\qquad \bU := \left(\bu_1\, | \cdots | \, \bu_M\right)\,, \label{xcompat}
\end{eqnarray}
\end{subequations}
where $\bb,\balpha \in \Real^{\cN}$ are given s.t. $\bU^{T} \bb = \bzero$, and $\bx \in \Real^{\cN}$ is the vector of unknowns. (\ref{linsys}) has a unique solution \cite[pp.~185--186]{trot}. Our goal is to develop a numerical iterative solver of (\ref{linsys}) that {\it requires $O(m)$ storage and $O(m \loget)$ operations to generate an $\ep$-accurate solution, for graphs arising in real-world applications.} The hidden constants should be small (in the hundreds, not millions). The solver should be parallelizable, and require a smaller cost to re-solve the system for multiple $\bb$'s -- a useful feature for time-dependent and other applications.

\subsection{Applications}
The linear system (\ref{linsys}) is fundamental to many applications, briefly sketched below. See Spielman's review \cite[\S 2]{icm10} for more details.

{\it Elliptic Partial Differential Equations (PDEs)} discretized on unstructured grids by the Finite Element Method (FEM) \cite{boman}, typically solved at each time step of a time-dependent Computational Fluid Dynamics (CFD) simulation \cite{fischer}. $\bx$ represents pressure, $w$ is the material's diffusion coefficient, and $\bb$ is a forcing term.

{\it Network flows.} The Maximum Flow and Minimum Cost Flow problems are linear programming problems that, if solved by interior point algorithms, reduce to solving a sequence of restricted Laplacian systems \cite{DS08, FG07}.

{\it Electrical networks.} In the study of electrical flow through a resistor network $G$, nodes are electrical components, $\bx$ is the electrical potential, $w_{uv}$ is the conductance between $u$ and $v$, and $\bb$ is a vector of external currents.

{\it A stepping stone toward the eigenproblem.} Our multilevel methodology is applicable to finding the smallest eigenpairs of $\bA$ with minor adaptations (cf. \S \ref{eigenvalue}). The Laplacian eigenproblem is central to graph regression and classification in machine learning \cite{hongkong_class, zhu_learning}; spectral clustering of images and graph embedding \cite{radu}; dimension reduction for genetic ancestry discovery \cite{ancestry}; and spectral graph theory. Of particular interest is the Fiedler value -- the smallest non-zero eigenvalue and the corresponding Fiedler eigenvector, which measure the algebraic connectivity of $G$ \cite[\S 1.1]{chung} and related to minimum cuts \cite{ding}. Although we believe it preferable to develop multiscale strategies for the original formulations of these graph problems, as demonstrated by the work \cite{alg_distance} (more examples are discussed in the papers \cite{class04, clustering06, nature}), a fast black-box eigensolver is a practical alternative.

\subsection{Related Work}
There are two main approaches to solving (\ref{linsys}): direct, leading to an exact solution; and iterative, which produces successive approximations $\tbx$ to $\bx$ and typically requires $O(\loget)$ iterations to achieve $\ep$-accuracy, namely,
\be
	\|\bx-\tbx\|_{\bA} \leq \ep \|\bx\|_{\bA}\,,\qquad \|\bx\|_{\bA} := \sqrt{E(\bx)}\,.
	\label{epaccuracy}
\ee

\subsubsection{Direct Methods}
The Cholesky factorization with a clever elimination order can be applied to $\bA$. A permutation matrix $\bP$ is chosen and factorization $\bP^{T} \bA \bP = \bL \bL^{T}$ constructed so that the lower triangular $\bL$ is as sparse as possible, using Minimum or Approximate Minimum Degree Ordering \cite{md, add} methods. Except for simple graphs, direct algorithms do not scale, requiring $O(n^{1.5})$ operations for planar graphs and $O(n^3)$ in general. Alternatively, fast matrix inversion can be performed in $O(n^{2.376})$ or combined with Cholesky, yet yields similar complexities \cite[\S 3.1]{icm10}.

\subsubsection{Iterative Methods: Graph Theoretic}
These are variants of the Preconditioned Conjugate (PCG) method \cite[\S 10.3]{gvl} that achieve (\ref{epaccuracy}) in $O(\sqrt{\kappa(\bA \bB^{-1})}$ $\loget)$ iterations for a preconditioner $\bB$; $\kappa$ is the finite condition number \cite[\S 3.3]{icm10}.

Recent works \cite{koutis, st06} have been focusing on multilevel graph-sparsifying preconditioners. The graph is repeatedly partitioned into sets of high-conductance nodes without removing too many edges. The complexity is a near-linear $O(m \log^2 n \log(1/\ep))$, guaranteed for any symmetric diagonally-dominant $\bA$. Unfortunately, no implementation are available yet, nor is there a guarantee on the size of the hidden constant.

\subsubsection{Iterative Methods: AMG}
Algebraic Multigrid (AMG) is a class of high-performance linear solvers, originating in the early 1980s \cite[\S 1.1]{guide}, \cite{geodetic}, \cite{rs86} and still under active development. During a setup phase, AMG recursively automatically constructs a multi-level hierarchy of increasingly coarser graphs by examining matrix entries, without relying on geometric information. The solve phase consists of standard multigrid cycles. AMG can be employed either as a stand-alone solver or as a PCG preconditioner \cite[App.~A]{trot}. Open-source parallel implementations such as Hypre \cite{hypre} and Trilinos-PETSc \cite{trilinos} are available. In classical AMG, the coarse set is a subset of $\cN$; popular alternatives are aggregation AMG \cite[App.~A.9]{trot},\cite{braess},\cite{blatt} and smoothed aggregation \cite{sa}, where the coarse nodes are aggregates of fine nodes.

AMG mainly targets discretized elliptic PDEs on unstructured grids \cite{fischer}. Recent works have been focusing on improving the coarsening and interpolation to increase the solution efficiency. These include Bootstrap AMG \cite[\S 17.2]{yes}, \cite{bamg} adaptive smoothed aggregation \cite{asa} and interpolation energy minimization schemes such as Olson and his associates' \cite{olson}. While these methods serve to increase AMG's scope and approach linear scaling for more systems, they are not designed for general graphs, as will be explained in \S \ref{considerations}. The present work aims at generalizing AMG to graph Laplacians and addresses peculiarities not encountered in traditional AMG applications.

\subsection{Our Contribution}
We present Lean Algebraic Multigrid (LAMG): a practical graph Laplacian solver. LAMG attains optimal efficiency: its setup phase requires $O(m)$ time and storage, and solve phase requires $O(m \log(1/\ep))$ operations per Right-Hand Side (RHS). An unoptimized \matlab LAMG implementation scaled linearly for \numgraphs real-world graphs with up to \maxedges edges, ranging from computational fluid dynamics to social networks. While we do not prove nor claim that LAMG works for every graph, these results will hopefully support its practical use.

LAMG is an aggregation-AMG algorithm \cite[App.~A.9]{trot} composed of lean components that significantly decrease setup time and memory usage and boost the solve phase efficiency. The key design decision is the choice of a {\it caliber-1} (piecewise-constant) interpolation between levels. Fast asymptotic convergence is achieved by (a) a new {\it relaxation-based node proximity} definition, which guides the aggregation and improves upon the algebraic distance defined by Ron et al. \cite{alg_distance}; and (b) an {\it energy correction} applied to coarse-level systems. We offer two alternatives: a flat correction similar to Braess' work \cite{braess}, yet resulting in a superior efficiency; or an adaptive correction via multilevel iterant recombination \cite[\S 7.8.2]{trot} that is even more efficient.

Importantly, the developed multilevel methodology is extensible to eigenvalue problems, other linear systems and other graph computational problems.

\subsection{Paper Organization}
We begin by elucidating the general issues of developing AMG for the graph Laplacian in \S \ref{considerations}. The LAMG algorithm is explained in \S \ref{algorithm}, comprising of a setup phase (\S \ref{setup}) and a solve phase (\S \ref{solve}). Numerical results are presented in \S \ref{results}. Future enhancements and extensions are outlined in \S \ref{extensions}.

\section{General Considerations}
\label{considerations}
It is important to first understand the pitfalls of existing AMG algorithms in general graphs and their remedies in LAMG.

\subsection{Node Proximity}
The construction of an effective coarse node set hinges upon defining which nodes in $\cN$ are ``proximal'', i.e., nodes whose values are strongly coupled in all smooth (low-energy) vectors \cite[p.~473]{trot}. Table~\ref{prox_def} lists three definitions.
\begin{table}[htbp]
\centering
\begin{tabular}{|l|l|}
    \hline
    Classical AMG & $|w_{uv}| / \max\left\{ \max_s |w_{us}|, \max_s |w_{sv}| \right\}$ \\ \hline
    1/Algebraic Distance & $ 1/\left(\max_{k=1}^K \left| x_u^{(k)} - x_v^{(k)} \right|\right)$ \\ \hline
    \multirow{3}{*}{Affinity} & $c_{uv} / \max\left\{ \max_{s \not = u} c_{us}, \max_{s \not = u} c_{sv}\right\}\,,$ \\ & $c_{uv} := \left|\left(X_u,X_v \right)\right|^2/\left(\left(X_u,X_u \right)^2 \left(X_v,X_v \right)^2\right)\,,$ \\
    & $(X_u,X_v):= \sum_{k=1}^K x^{(k)}_u x^{(k)}_v$ \\
    \hline
\end{tabular}
\caption{Comparison of node proximity measures. Nodes are defined as ``close'' when the measure exceeds a certain threshold.}
\label{prox_def}
\end{table}
\begin{figure}[htbp]
    \centering
    \begin{center}
    \begin{tabular}{cc}
      \includegraphics[height=1.2in]{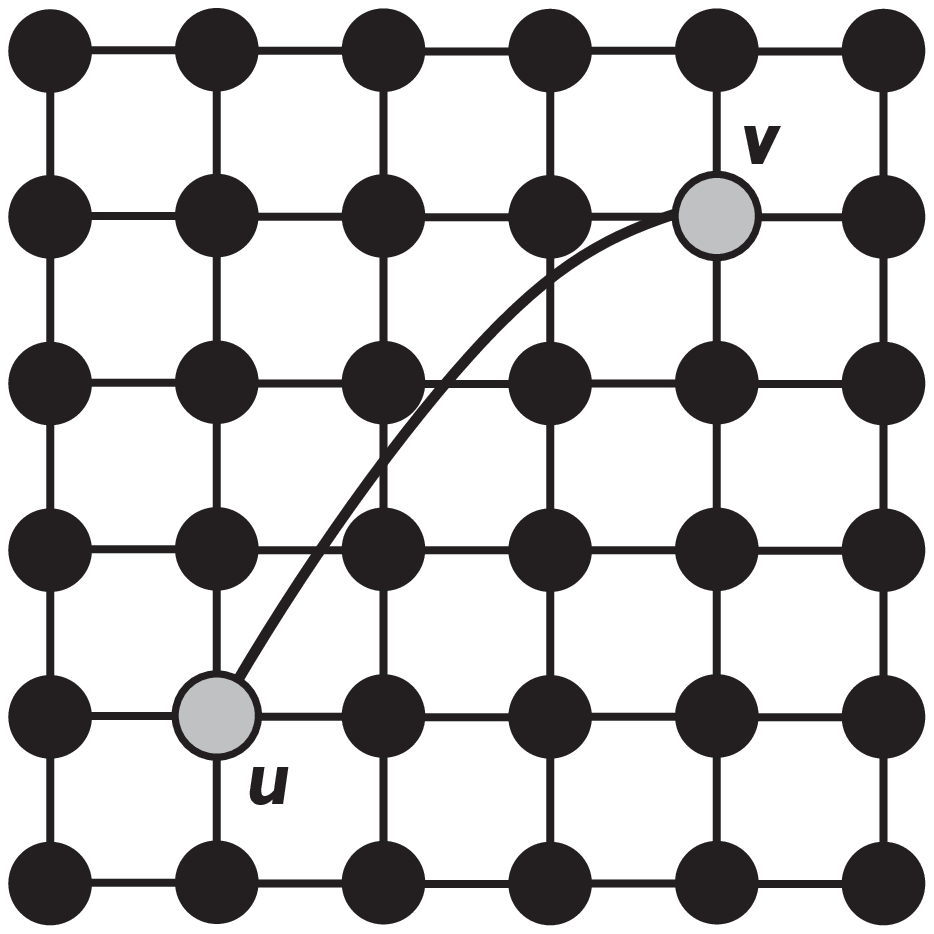}
      &
      \includegraphics[height=1.15in]{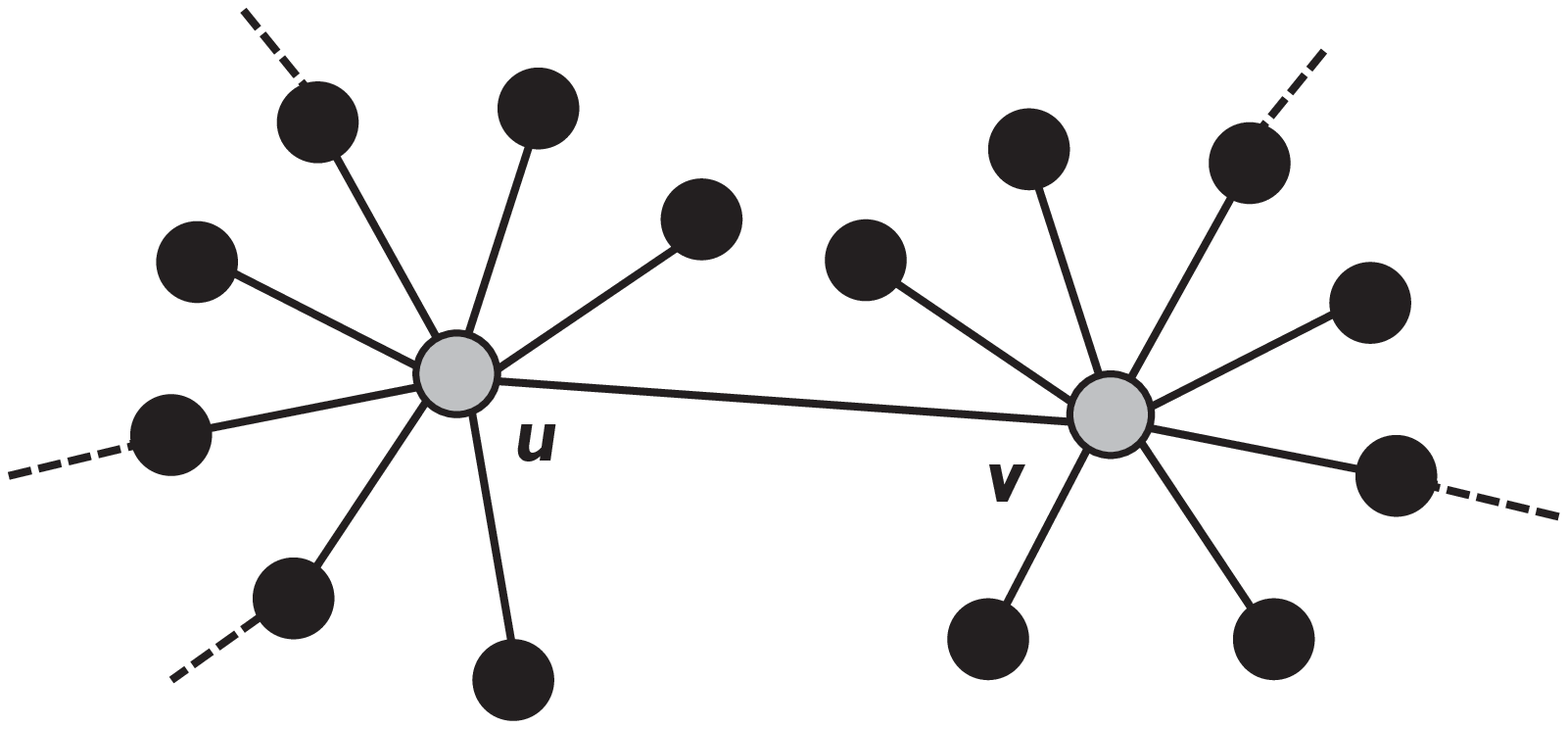}
      \\
      (a) & (b) \\
    \end{tabular}
    \end{center}
    \caption{Graph instances that present aggregation difficulties. All graph edge weights are $1$. (a) A 2-D grid with an extra link. (b) Two connected suns. A sun is a high-degree node.}
    \label{two_examples}
\end{figure}

Classical AMG defines proximity based on matrix entries (Table~\ref{two_examples}, top row). While this has worked well for coarsening discretized scalar elliptic PDEs, it leads to wrong aggregation decisions in non-local graphs. In a grid graph with an extra link between distant nodes $u$, $v$ (Fig.~\ref{two_examples}a), $u$ and $v$ become proximal and may be aggregated. Unless $w_{uv}$ is outstandingly large, this is undesirable because $u$ and $v$ belong to unrelated milieus of the grid.

This problem is overcome by the algebraic distance measure introduced by Ron et al. \cite{alg_distance} (Table~\ref{prox_def}, middle row; a related definition is used in the work \cite{alg_distance_anis}). A set of $K$ relaxed Test Vectors (TVs) $\bx^{(1)},\dots,\bx^{(K)}$ is generated, each obtained by applying several Gauss-Seidel (GS) relaxation sweeps to $\bA \bx = \bzero$, starting from $\mathrm{random}[-1,1]$ and normalizing the result. Yet this definition falls pray to a graph containing two connected {\it suns} (high-degree nodes) $u$ and $v$, each of which connecting many satellites (Fig.~\ref{two_examples}b). For each $k$, the value $x^{(k)}_u$ is an average over a large neighborhood of nodes whose size increases with the number of sweeps, and will be small. Similarly, $x^{(k)}_v$ is small, so every pair of suns is always proximal even though they may represent distant node clusters.

LAMG's proximity measure is the {\it affinity} (Table~\ref{prox_def}, last row), which also relies on the same TVs, yet is scale-invariant and correctly assesses both cases of Fig.~\ref{two_examples} as well as many other constellations.

\subsection{Interpolation Caliber}
Textbook multigrid convergence for the Poisson equation requires that the interpolation of corrections $\bP$ be second-order \cite[\S 3.3]{mg_theory}. The analogous AMG theory implies a similar condition on the interpolation accuracy of low-energy errors. While a piecewise-constant $\bP$ is acceptable in a two-level cycle, it is insufficient in V-cycles; W-cycles are faster but more costly \cite[p.~471]{trot}.

Constructing a second-order $\bP$ is a challenge even for grid graphs. Often a first-order interpolation is constructed, followed by a smoothing step \cite{sa} or a scheme that reduces the interpolation's energy \cite{bamg, olson}. The situation is exacerbated in graphs with no geometric structure: the interpolation order is undefinable. Even had the graph's effective dimension $d$ been known, the required interpolation {\it caliber}, i.e., the number of coarse nodes used to interpolate a fine node, would grow with $d$ and result in unbounded interpolation complexity. Finally, choosing a proper interpolation set (whose ``convex hull'' contains the fine node) is a complex, costly endeavor \cite{bamg}.

In contrast, LAMG employs a caliber-1 (piecewise-constant) $\bP$, and corrects the energy of the coarse-level Galerkin {\it operator} to maintain good convergence (in practice, the correction is actually applied to the coarse RHS). This could not have been achieved within the variational setting, which only permits modifying $\bP$. Whereas high-caliber works focus on optimizing $\bP$, here the barrier to fast convergence is the coarse-to-fine operator energy ratio. Our contribution is an algorithm that yields a small energy ratio, which translates into optimal efficiency; cf. \S\S \ref{aggregation}--\ref{energy_correction}. The Compatible Relaxation (CR) performance predictor \cite{cr_etna, cr_oren, cr_james}, \cite[\S\S 14.2--14.3]{guide} is not relevant for low interpolation accuracy; the energy ratio is a better predictor.

\subsection{Coarse-level Fill-in}
Frequently, the coarse-level matrices in AMG hierarchies become increasingly dense. This is a result of a poor aggregation, a high-caliber $\bP$, or both: many fine nodes whose neighbor sets are disjoint are aggregated, creating additional edges among coarse-level aggregates. This renders the ideally-accurate interpolation irrelevant, because the actual cycle {\it efficiency} (error reduction per unit work) is small even though convergence may be rapid. While fill-in is often manageable in grid graphs because the coarse graphs are still local, it is detrimental in non-local graphs.

LAMG's interpolation is designed to create an insignificant fill-in: the sparser $\bP$, the sparser the Galerkin operator $\bP^{T} \bA \bP$. The affinity-based aggregation (\S \ref{aggregation}) also helps, as it tends to aggregate nodes with many common neighbors. The cycle work is further controlled by a fractional cycle index \cite[\S 6.2]{guide} between 1 and 2; cf. \S \S \ref{setup},\ref{solve}.

Occasionally, the interpolation caliber may be slightly increased as long as the number of coarse edges does not become too large; see \S \ref{improvements}.

\subsection{Extenuating Circumstances}
\label{exten}
Specific properties of the graph Laplacian can be exploited to simplify the LAMG construction.
\bi
    \item GS is an effective smoother in SPS systems \cite[\S 1]{guide}. Also, the quadratic form (\ref{energy}) is handy for discerning and solving the energy inflation problem (\S \ref{energy_correction}).
    \item Since $\bA$ has zero row sums, its null-space eigenvectors are constant in every connected component. Under the general AMG assumption that {\it all near}-null-space errors can be fitted by a {\it single} interpolation from a coarse level \cite[p.~8]{guide}, the caliber-1 interpolation weights are apriori set to $1$. (This assumption is easily verified for Laplacians with bounded node degrees \cite[p.~439]{trot}. It is violated in wave equations, where multiple coarse grids are required \cite{wave_accuracy}.)
    \item Some graph locales are effectively one-dimensional: many nodes have degree 2-3. Such nodes can be quickly eliminated similarly to the paper \cite{koutis} (\S \ref{elimination}).
    \item In other graphs GS is an efficient solver and no coarsening is required. These include complete graphs, star graphs and expander graphs \cite[\S 1]{icm10}.
\ei

\subsection{Data Structures}
The efficiency of iterative methods for (\ref{linsys}) depends on a proper storage format for $\bA$, as they compute many matrix-vector products. AMG algorithms also require other types of operations, e.g., removing a row and a corresponding column upon node aggregation. Standard choices such as Compressed Column storage (CCS) \cite{templates, matlab_sparse} may not be optimal for these operations, especially in graphs containing suns. In our implementation, for instance, we do not recompute affinities upon a TV update in Algorithm~\ref{aggregationStage} because that would require a slow \matlab CCS matrix row update.

Optimized architectures for such operations that best utilize the available hardware (e.g., CPU or GPU) may be pursued \cite{fischer}.

\section{The LAMG Algorithm}
\label{algorithm}
Veracious to the algebraic multigrid framework, LAMG recursively constructs a hierarchy of $L$ increasingly coarser Laplacian systems (``levels'') $\bA^l \bx^l=\bb^l$, $l=1,\dots,L$ during a setup phase. The finest level is the original system, $\bA^1:=\bA, \bb^1:=\bb$. The setup phase depends on $\bA$ only, and constructs a sequence $\{(\bA^l,\bP^l)\}_{l=2}^L$, where $\bA^l$ is an $n_l \times n_l$ Laplacian and $\bP^l$ is an $n_{l-1} \times n_l$ interpolation matrix from level $l$ to $l-1$. We also use the notation $G^l=(\cN^l,\cE^l)$ for the graph corresponding to $\bA^l$, although it is not explicitly stored by our code. Each level $l \geq 2$ is of type \elim (obtained from the next-finer level by exact node elimination) or \agg (a caliber-1 node aggregation of fine-level nodes). The vectors $\{\bb^l\}_{l=2}^L$ are computed during the solve phase.

\subsection{Setup Phase}
\label{setup}
\begin{figure}[htbp]
    \centering
    \begin{center}
      \includegraphics[height=3.5in]{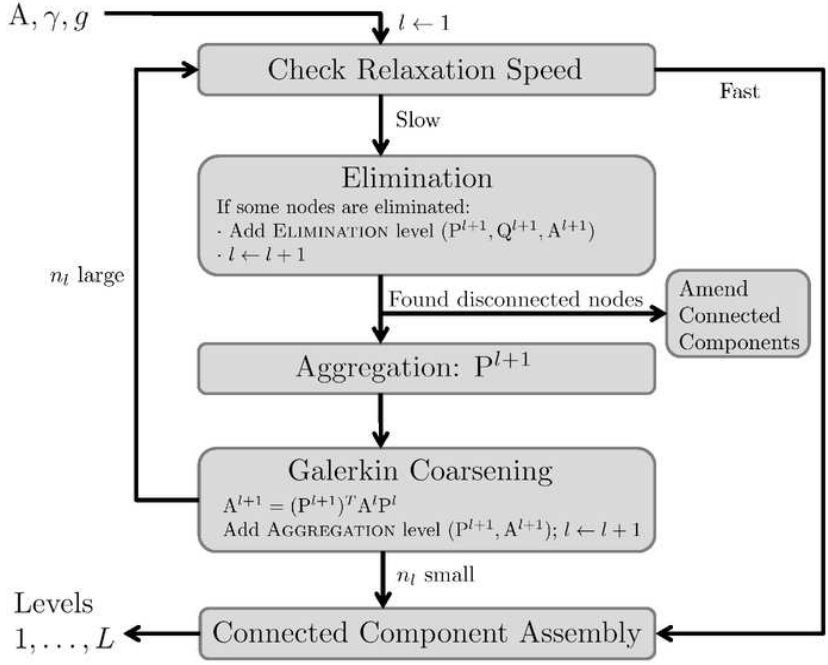}
    \end{center}
    \caption{LAMG setup phase flowchart.}
    \label{flowchart}
\end{figure}

The setup flow is depicted in Fig.~\ref{flowchart}. It requires two inputs:
\bi
    \item {\it Cycle index $\gamma \geq 1$} to be employed at most levels of subsequent solution cycles.
    \item {\it Guard $0 < g < 1$}, which bounds coarsening ratios and controls the cycle work.
\ei
In our program, $\gamma=1.5$ and $g=.7$; these values are discussed in \S \ref{solve}.

Given the current coarsest Laplacian $A^l$, we first estimate the speed of relaxation (\S \ref{relax}), and terminate if it is fast enough. Otherwise, disconnected and low-degree nodes are eliminated from $G^l$ (\S \ref{elimination}). If such nodes are found, an \elim level is added to the hierarchy, and $G^l$'s disconnected nodes are amended to the list of $G$'s connected graph components. Let $l$ be the new coarsest level; next, $\cN^l$ is aggregated to form $\cN^{l+1},P^{l+1}$ (\S \ref{aggregation}). Finally, an energy-corrected Galerkin operator $A^{l+1}$ is constructed (\S \ref{energy_correction}). Coarsening is repeated until $n_l$ is small enough or until relaxation is fast. Finally, $G$'s connected components are assembled (\S \ref{connected_comp}).

\subsection{Relaxation}
\label{relax}
Let $l$ be the currently processed level during setup. For simplicity, we omit the $l$-superscripts in this section. A Gauss-Seidel (GS) relaxation sweep for $\bA \bx=\bb$ is defined by the successive updates \cite[\S 1.1]{guide}
\be
    \text{For } u=1,\dots,n\,,\qquad x_u \leftarrow \frac{b_u - \sum_{v \in \cA_u} a_{uv} x_v}{a_{uu}}\,.
    \label{gs}
\ee
A GS solve iteration is a GS sweep followed by subtracting the mean of $\bx$ from all $\{x_u\}_u$. Starting from $\bx_0 = \mathrm{random}[-1,1]$, we apply $\nu=15$ GS solve iterations to $\bA \bx = \bzero$ and estimate the Asymptotic Convergence Factor (ACF) by $\rho:=\|\bx_{\nu}\|_2/\|\bx_{\nu-1}\|_2$, where $\bx_i$ is the $\kth{i}$ iterant and $\|\cdot\|_2$ is the $l_2$ norm. If $\rho \leq .7$, $l$ becomes the final coarsest level $L$. A slower GS solve implies that the graph is ``stiff'', i.e. there exist low-energy errors whose magnitude is not reflected by their residuals \cite[\S 1.1]{guide}, hence further coarsening is necessary.


The speed check is skipped if $G$ has disconnected ($0$-degree) nodes, since a division by $0$ would occur in (\ref{gs}).

\subsection{Low-degree Node Elimination}
\label{elimination}

First, we eliminate from $G$ all disconnected nodes $\cZ$, and a set of low-degree nodes $\cF$. $\cZ$-elimination is mandatory for GS to be properly defined on the remaining graph. The $\cF$-elimination ventures to reduce $n$ while not significantly increasing $m$, at a small cost. This removes the 1-D part of the graph and enhances the efficiency of subsequent \agg levels (cf. \S \ref{energy_correction}).

\subsubsection{Choosing $\cF$}
$\cF$ is an independent set of nodes $u$ with degree $|\cA_u| \leq 4$. Eliminating a node connects all its neighbors; hence $\cF$-nodes with $|\cA_u| \leq 3$ do not increase $m$ (Table~\ref{elim_cases}). When $|\cA_u| = 4$, $m$ might be increased by at most $2$. However, we assume that this is unlikely to happen for many $\cF$ nodes, hence eliminate those as well (a costly alternative is to monitor the future change in $m$ and only eliminate nodes that do not increase it). Larger $|\cA_u|$ values result in an impractical fill-in.

\begin{table}[htbp]
\footnotesize
\centering
\begin{tabular}{|c|c|l|l|}
  \hline
  \multirow{2}{*}{$|\cA_u|$} &
  \multirow{2}{*}{$\Delta_m$} &
  Before & After \\
  && Elimination & Elimination \\ \hline
  $1$ & $-1$ &
  \begin{minipage}{.15\linewidth}
     \includegraphics[height=.15in]{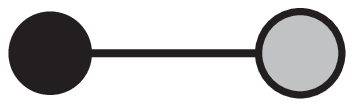}
  \end{minipage}
  &
  \begin{minipage}{.15\linewidth}
     \includegraphics[height=.17in]{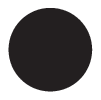}
  \end{minipage}
  \\ \hline
  $2$ & $-1$ &
  \begin{minipage}{.15\linewidth}
     \includegraphics[height=.18in]{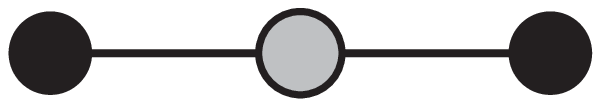}
  \end{minipage}
  &
  \begin{minipage}{.15\linewidth}
     \includegraphics[height=.12in]{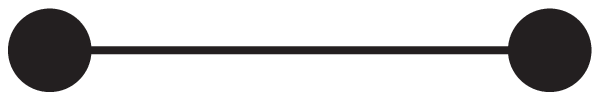}
  \end{minipage}
  \\ \hline
  $3$ &  $0$ &
  \begin{minipage}{.15\linewidth}
     \includegraphics[height=.38in]{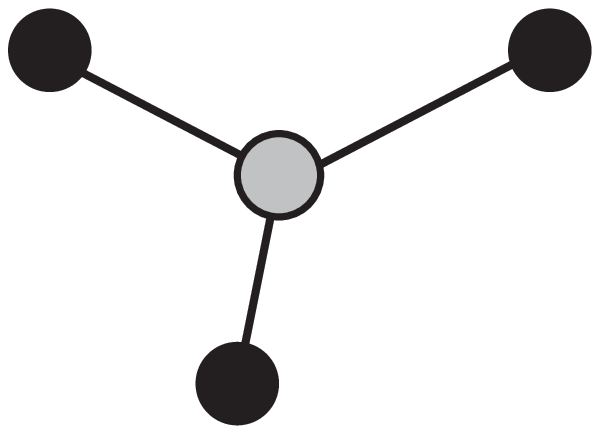}
  \end{minipage}
  &
  \begin{minipage}{.15\linewidth}
     \includegraphics[height=.38in]{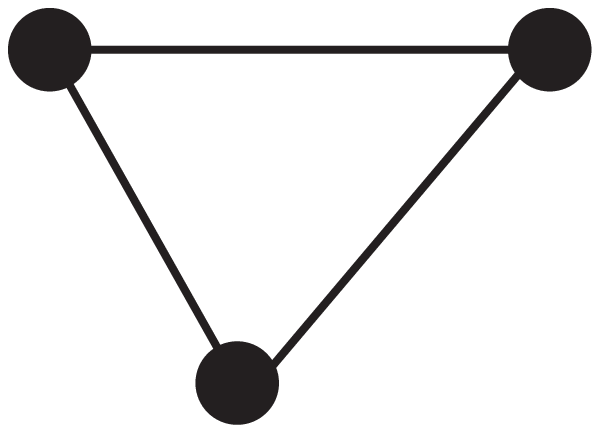}
  \end{minipage}
  \\ \hline
  $4$ & $+2$ &
  \begin{minipage}{.15\linewidth}
     \includegraphics[height=.38in]{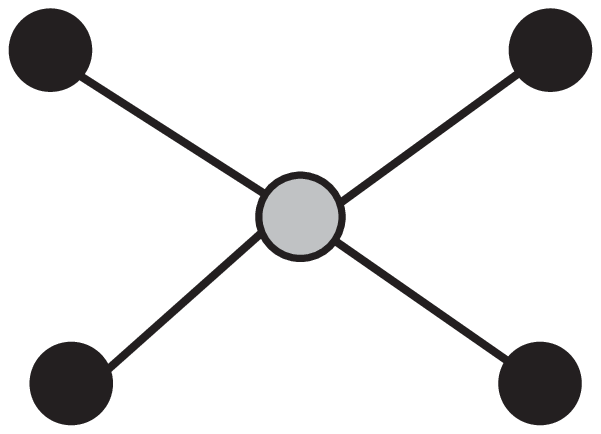}
  \end{minipage}
  &
  \begin{minipage}{.15\linewidth}
     \includegraphics[height=.38in]{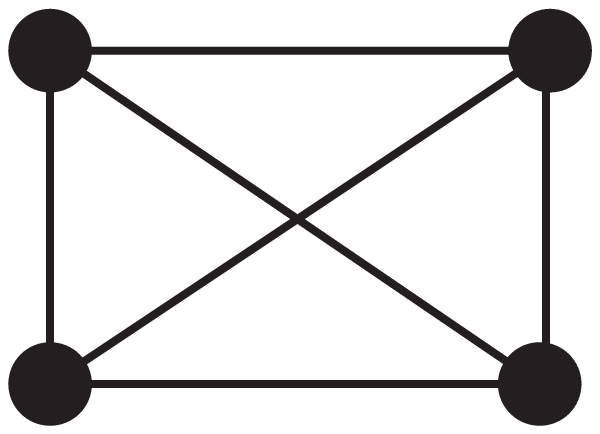}
  \end{minipage}
  \\ \hline
\end{tabular}
\caption{Types of nodes $u \in \cF$. For each type, the sub-graph of $u$ (gray circle) and $\cA_u$ (black circles) is depicted before and after elimination in the worst fill-in case. $\Delta_m$ is the maximum increase in $m$.}
\label{elim_cases}
\end{table}

The pseudo-code of choosing $\cF$ is given in Algorithm~\ref{low_degree}. Nodes are sequentially visited; when an $\cF$ node is marked, its neighbors become ineligible for inclusion in $\cF$, which guarantees the set's independence.
\begin{algorithm}
\caption{$\cF = LowDegreeNodes(\bA)$}
\label{low_degree}
\begin{algorithmic}[1]
    \STATE $\cN \leftarrow \nodes(\bA)$, $n \leftarrow \sz(\bA)$, $\cU \leftarrow \left\{ u \in \cN: 1 \leq |\cA_u| \leq 4 \right\}.$
    \STATE For each $u \in \cU$, $\visited(u) \leftarrow \notvisited$.
    \FOR{each $u \in \cU$}
        \IF{($\visited(u) = \notvisited$)}
            \IF[$u$ may be eliminated]{$\not \exists v \in \cA_u: \visited(v) = \fnode$}
                \STATE $\visited(u) \leftarrow \fnode$.
                \STATE For each $v \in \cA_u$, $\visited(v) \leftarrow \notelim$.
            \ELSE[$u$ has an $\cF$-neighbor]
                \STATE $\visited(u) \leftarrow \notelim$.
            \ENDIF
        \ENDIF
    \ENDFOR
    \STATE return $\left\{ u \in \cU: \visited(u) = \fnode \right\}$.
\end{algorithmic}
\end{algorithm}

\subsubsection{Elimination Equations}
Let $\cC := \cN \ (\cZ \cup \cF)$ denote the rest of the nodes, and permute rows and columns so that (\ref{system}) becomes
\be
    \begin{pmatrix}
    \bzero & \bzero & \bzero \\
    \bzero & \bA_{\cF\cF} & \bA_{\cF\cC} \\
    \bzero & \bA_{\cF\cC}^{T} & \bA_{\cC\cC}
    \end{pmatrix}
    \begin{pmatrix}
    \bx_{\cZ} \\
    \bx_{\cF} \\
    \bx_{\cC} \\
    \end{pmatrix}
    =
    \begin{pmatrix}
    \bb_{\cZ} \\
    \bb_{\cF} \\
    \bb_{\cC} \\
    \end{pmatrix}
    \,.
    \label{block}
\ee
Note that $A_{\cF\cF}$ is diagonal. Block Gaussian elimination of the $\cZ$- and $\cF$-blocks transforms (\ref{block}) into
\begin{subequations} \label{reduced}
\begin{eqnarray}
    && \bx_{\cZ} = \bzero \label{xz} \\
    && \bx_{\cF} = \bA_{\cF\cF}^{-1} \left( \bb_{\cF} - \bA_{\cF\cC} \bx_c \right) \label{xf} \\
    && \left( \bA_{\cC\cC} - \bA_{\cF\cC}^{T} \bA_{\cF\cF}^{-1} \bA_{\cF\cC} \right) \bx_{\cC} = \bb_{\cC} - \bA_{\cF\cC}^{T} \bA_{\cF\cF}^{-1} \bb_{\cF} \label{xc}
\end{eqnarray}
\end{subequations}
(\ref{xc}) is the Schur complement system \cite{schur_book}, and (\ref{xz})--(\ref{xf}) recover the rest of the variables once it is solved. (\ref{reduced}) can be written more succinctly as
\begin{subequations} \label{elim}
\begin{eqnarray}
    && \bx = \bP \bx^c + \bQ \bb \label{elim_x} \\
    &&  \bP := \bPi \left( \bzero_{\cZ}\,, - \bA_{\cF\cC}^{T} \bA_{\cF\cF}^{-1}\,,\bI_{\cC} \right)^{T}\,,\quad \bQ := \bPi \text{\;diag}\left\{ \bzero_{\cZ}\,, \bA_{\cF\cF}^{-1} \,, \bzero_{\cZ} \right\}
     \label{elim_p} \\
    && \bA^c \bx^c = \bb^c\,, \qquad \bA^c := \bP^{T} \bA \bP\,,\quad \bb^c := \bP^{T} \bb\,,\quad \bx^c := \bx_{\cC}\,. \label{elim_coarse}
\end{eqnarray}
\end{subequations}
$\bPi$ is a permutation matrix such that $\bPi^{T} \bx$ lists the values of $\bx$ at $\cZ$ nodes, $\cF$ nodes and $\cC$ nodes in this order. (\ref{elim_coarse}) is a smaller Laplacian system for which another elimination stage is performed: new $\cZ$ and $\cF$ are identified, leading to a still-coarser (\ref{elim_coarse}), and so forth, until $\cZ=\emptyset$ and $\cF$ is relatively small (cf. Algorithm~\ref{elim_algorithm}; in principle, one can proceed until $\cF = \emptyset$ and maintain linear time by scanning only the neighbors of the previous stage's $\cF$-set in $LowDegreeNodes()$, but we haven't yet implemented this feature). For $q$ elimination stages, the operators $\{\bP_i,\bQ_i\}_{i=1}^q$ are lumped into composite $\bP=\bP_q \bP_{q-1} \cdot \ldots \cdot \bP_1$ and $\bQ$ to avoid storing a Laplacian $\bA^c$ per stage. For instance, when $q$ is large yet each stage eliminates only few nodes, storage would be prohibitively high.

\begin{algorithm}
\caption{$[\bA^c,\bP,\bQ,\{\bP_i,\bQ_i\}_{i=1}^q] = Elimination(\bA)$}
\label{elim_algorithm}
\begin{algorithmic}[1]
    \STATE $\bA^c \leftarrow \bA$, $\cN^c \leftarrow \nodes(\bA)$, $n_c \leftarrow |\cN^c|$, $\bQ \leftarrow \bzero_{n_c}$, $\bP \leftarrow \bI_{n_c}$, $q \leftarrow 0$.
    \WHILE{($n_c > 1$)}
        \STATE $\cZ \leftarrow \left\{ u \in \cN^c: |\cA^c_u| = 0 \right\}$, $\cF \leftarrow LowDegreeNodes(\bA^c)$.
        \IF{(($\cZ = \emptyset$) and ($|\cF| < .01 n_c$))}
            \STATE return $\bA^c,\bP,\bQ,\{\bP_i,\bQ_i\}_{i=1}^q$.
        \ENDIF
        \STATE $q \leftarrow q+1$, $\cC \leftarrow \cN \ (\cZ \cup \cF)$.
        \STATE $\bPi \leftarrow$ permutation matrix that transforms $\cN^c$ to the order $(\cZ,\cF,\cC)$.
        \STATE $\bQ_i \leftarrow \bPi \cdot \diag\{ \bzero_{\cZ}\,, (\bA^c_{\cF\cF})^{-1}\,, \bzero_{\cF}\}$,
        $\bR \leftarrow (\bzero_{\cZ}\,, \bzero_{\cF}\,, \bI_{\cC})$.
        \STATE $\bP_i \leftarrow \bPi \cdot (\bzero_{\cZ}\,, - (\bA^c_{\cF\cC})^{T} (\bA^c_{\cF\cF})^{-1}\,,\bI_{\cC})^{T}$.
        \STATE $\bQ \leftarrow \bQ + \bP \bQ_i \bR$, $\bP \leftarrow \bP_i \bP$,
        $\bA^c \leftarrow \bP^{T} \bA^c \bP$, $\cN^c \leftarrow \nodes(\bA^c)$, $n_c \leftarrow |\cN^c|$.
    \ENDWHILE
\end{algorithmic}
\end{algorithm}

If no nodes are eliminated during $Elimination()$, $l$ remains the coarsest level, otherwise $(\bA^c,\bP)$ defines the next level, $l+1$. Either way, we again denote the new coarsest system by $\bA \bx = \bb$ to be further coarsened in \S\S \ref{aggregation}--\ref{energy_correction}.

\subsection{Aggregation}
\label{aggregation}
\begin{figure}[htbp]
    \centering
    \begin{center}
    \begin{tabular}{cc}
      \includegraphics[height=1.2in]{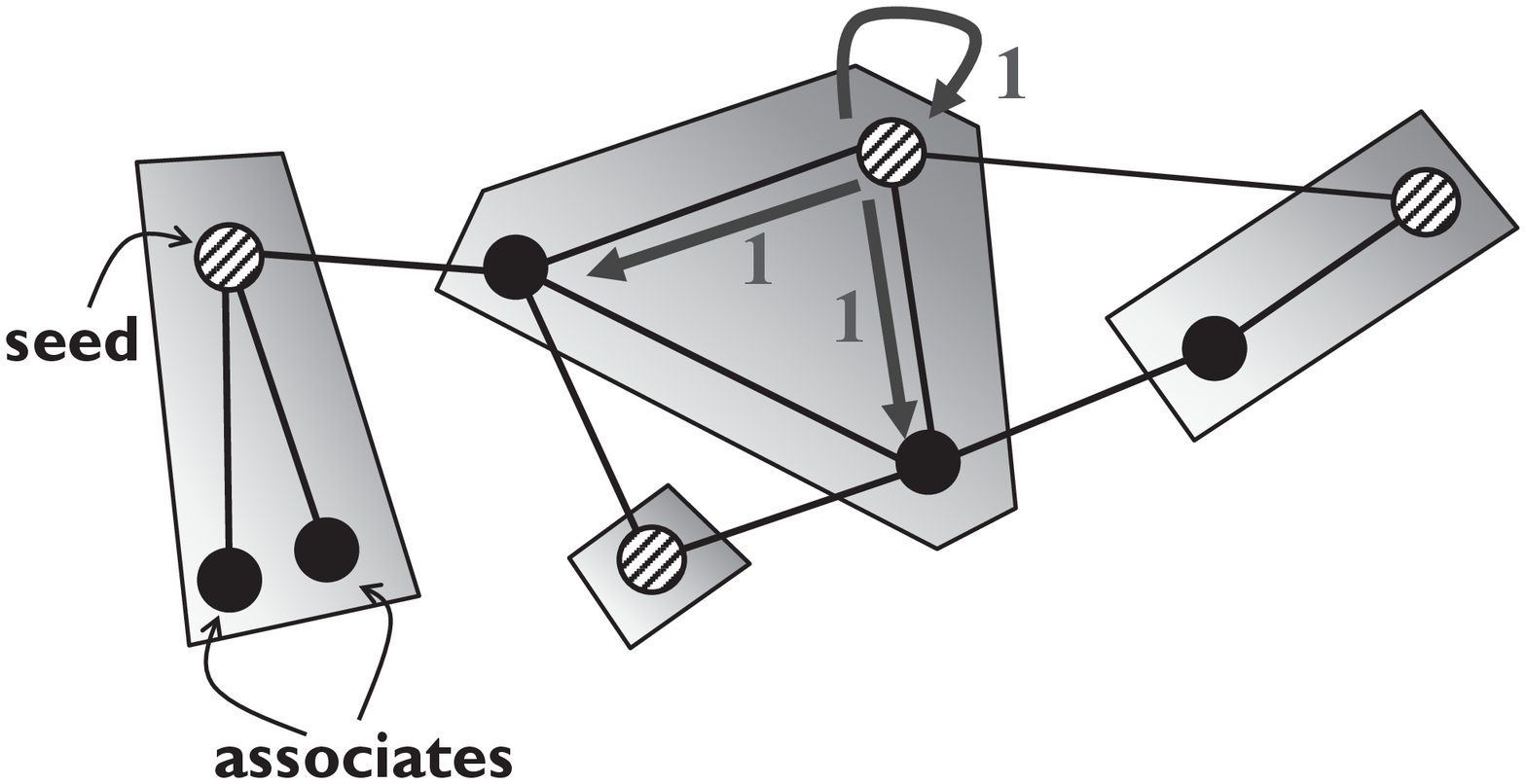}
      &
      \includegraphics[height=1.1in]{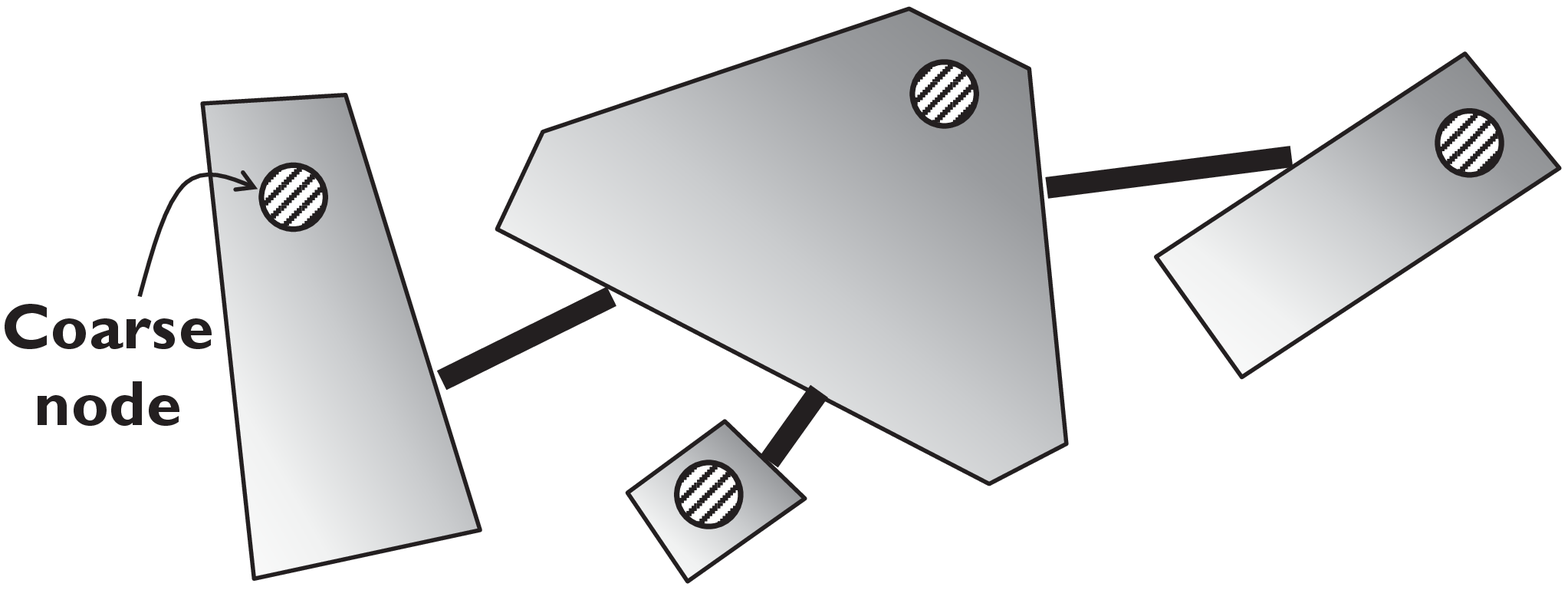}
      \\
      (a) & (b) \\
    \end{tabular}
    \end{center}
    \caption{(a) An aggregation. Seeds are the semi-filled nodes. The interpolation copies the coarse seed value to its associates. (b) The coarse graph. Edge weights are computed by (\ref{galerkin_entry}).}
    \label{agg}
\end{figure}

Eq.~(\ref{system}) is equivalent to the quadratic minimization
\be
    \bx = \mymin{\by} \Etot(\by)\,,\quad \Etot(\bx) := \frac12 \bx^{T} \bA \bx - \bx^{T} \bb\,.
    \label{quad}
\ee
Let $\tbx$ be an approximation to $\bx$ after several relaxation sweeps. The remaining error $\bee := \bx-\tbx$ of ``special nature'': its normalized residuals are much small than its magnitude \cite[\S 1.1]{guide} (assuming its mean of every connected component has been subtracted). Such errors are called {\it algebraically smooth}, or smooth for short, and are approximated by an interpolation from a coarse level, $\bee \approx \bP \bee^c$. The variational correction scheme \cite[\S 4.5]{guide} finds the optimal correction in the energy norm, namely,
\be
    \bee^c = \argmin{\by^c} \Etot\left(\tbx + \bP \by^c \right)\,;
    \label{min_ec}
\ee
The normal equations of this minimization are the Galerkin coarsening
\be
    \bA^c \bee^c = \bb^c\,,\quad \bA^c := \bP^{T} \bA \bP \,,\quad \bb^c := \bP^{T} \left(\bb - \bA \tbx\right) \,.
    \label{galerkin}
\ee
Once an approximation $\tbee^c$ to $\bee^c$ is computed, the fine-level approximation is corrected:
\be
    \tbx \leftarrow \tbx + \bP \tbee^c\,.
    \label{correction}
\ee
This algorithm depends only on $\bP$. Our particular $\bP$ has caliber $1$, which is equivalent to partitioning $\cN$ into $n_c$ non-overlapping {\it aggregates} $\{\cT_U\}_{U \in \cN^c}$, where $\cT_U$ is the set of $\bee_u$'s interpolated from $\bee_U$ \cite[App.~9]{trot} with unit weights (cf. \S \ref{exten}) and $\cN^c := \{1,\dots,n_c\}$. Each aggregate is composed of a {\it seed} node and zero or more {\it associate} nodes (Fig.~\ref{agg}). For simplicity, the seeds can be thought of as the coarse nodes, although $\bee_U$ is just a degree of freedom that can be interpreted differently. The Galerkin operator computation is simplified and involves additions only:
\be
   \bA^c_{UV} = \sum_{u \in \cT_U} \sum_{v \in \cT_V} a_{uv}\,,\qquad U,V \in \cN^c\,.
   \label{galerkin_entry}
\ee
Since $\bP$ has unit row sums, $\bA^c$ is SPS and zero-sum and hence a Laplacian.

Intuitively, nodes should be aggregated together if their values are ``close''. The next section introduces a node proximity measure.

\subsubsection{Affinity}
\label{affinity}
Insofar as coarsening concerns the space of smooth error vectors $\bx$, nodes $u$ and $v$ should be aggregated only if $\bx_u$ and $\bx_v$ are highly correlated for all such $\bx$. We generate $K$ Test Vectors (TVs) $\bx^{(1)},\dots,\bx^{(K)}$ -- a sample of this space \cite[\S 17.2]{yes}, \cite{alg_distance}. Each TV is the result of applying $\nu$ GS relaxation sweeps to $\bA \bx = \bzero$, starting from $\mathrm{random}[-1,1]$. Let $\bX_{n \times K} := \left(\bx^{(1)}\, | \cdots | \, \bx^{(K)}\right)$.

Since TVs are used to derive a coarsening of a modest coarsening ratio $n_c/n$ (typically $\frac13$ to $\frac12$), they need not be overly smooth nor numerous: $\nu=3$ is used at all levels; $K=8$ TVs are employed to form the finest aggregation level, and an extra TV is added at each aggregation level. This incurs a small additional cost and seems useful, as coarse-level graphs are often increasingly more dense and complex. (The asymptotic GS solve vector obtained in \S \ref{relax} can be reused as a TV to save work.)

The {\it affinity} $c_{uv}$ between nodes $u$ and $v$ is the goodness of fit of fitting the linear model $x_v \approx x_u$ to TV values:
\be
    c_{uv} := \frac{\left|\left(X_u, X_v\right)\right|^2}{\left(X_u, X_u\right)^2 \left(X_v, X_v\right)^2}\,,
    \qquad \left( X,Y \right) := \sum_{k=1}^K X^{(k)} Y^{(k)}\,.
    \label{cuv}
\ee
For every $u,v$, $c_{uu} = 1$, $0 \leq c_{uv} \leq 1$ and $c_{uv}=c_{vu}$. The affinity measures the {\it strength of connection}: the larger $c_{uv}$, the closer $u$ and $v$. Let $0 \leq \delta \leq 1$; nodes $u$ and $v$ are called {\it $\delta$-affinitives} if
\be
    c_{uv} \geq \delta \max\left\{ \max_{s \not = u} c_{us}\,, \max_{s \not = v} c_{sv}\right\}\,.
    \label{delta_aff}
\ee
This is analogous to the classical AMG definition (Table~\ref{prox_def}, top row) and works well in practice. Notwithstanding, alternative definitions that account for the closeness of $c_{uv}$ to $1$ need to be explored.


Drawing an analogy to statistics, $x_u$ can be thought of as a random variable; $c_{uv}$ is the coefficient of determination ($R^2$) of linearly regressing $x_u$ on $x_v$ \cite{regression} using the TV sample. This interpretation leads to several observations.
\bi
    \item $1-c_{uv}$ is an alternative definition of the algebraic distance between $u$ and $v$ \cite{alg_distance}. It is related to {\it geometric distance} in grid graphs. For instance, in the 1-D discretized Poisson equation with meshsize $h=1/n$ whose nodes are located at $\{u h\}_{u=1}^n$, $1-c_{uv} \propto ((u-v)h)^2$ in the limit of $h \rightarrow 0$.
    \item $c_{uv}$ is invariant to scaling $X_u$ and $X_v$, which is vital in the two suns case (Fig.~\ref{two_examples}b).
    \item $c_{uv}$ is unbiased by the sample size $K$, because we use the exact means ($0$) of $X_u$ and $X_v$ over all error vectors, rather than the sample means $\bar{x}_u := \sum_{k=1}^K x^{(k)}_u$ and $\bar{x}_v$. Indeed, the initial TVs are uniformly distributed, so the probabilities of starting from $\bx$ and from $-\bx$ are equal; since relaxation is a linear process, the probability of encountering a relaxed TV $\bx$ equals the probability of $-\bx$, thus the mean of $x_u$ over all possible TVs is $0$.
    \item In principle, Least-squares regression weights should be proportional to each measurement's reciprocal variance \cite{regression}. The variance of a TV is proportional to an appropriate norm of its normalized residuals, so that smoother the TV, the larger its weight \cite[\S 17.2]{yes}, \cite{bamg}. In our case, all vectors have equal weights because they have the same smoothness level, i.e., comparable normalized residuals. This saves work and avoids biases that may occur if an improper weighting is applied. (For example, no TV should be given a much larger weight than all others, in which case it would dominate the regression and lead to a nonsensical $c_{uv}$.)
\ei
The work \cite{alg_distance} also defines the affinity of node $u$ to a node {\it set} $\cV$. Here, the analogue is the coefficient of multiple determination of regressing $x_u$ on $\{x_v\}_{v \in \cV}$,
\be
    c_{u\cV} := \frac{(\hat{X}_u,\hat{X}_u)}{(X_u,X_u)}\,,\quad
    \hat{X}_u := \sum_{v \in \cV} \hat{q}_v X_v\,,\quad
    \hat{\bq} := \argmin{\bq} \left\|X_u-\sum_{v \in \cV} q_v X_v\right\|^2\,.
    \label{affinity_set}
\ee
As a byproduct we obtain the interpolation coefficients $\hat{\bq}$ from $\cV$ to $u$. (\ref{affinity_set}) reduces to (\ref{cuv}) for $\cV=\{v\}$, where the corresponding regression model is $x_u \approx \hat{q} x_v$ with $\hat{q} = (X_u,X_v)/(X_v,X_v)$. Thus (\ref{cuv}) also works for non-zero row-sum M-matrices, e.g., restricted Laplacians \cite[\S 2]{icm10}. In the Laplacian, $\hat{q}$ is abandoned in favor of the theoretically known interpolation weight $1$ (cf. \S \ref{exten}); in non-zero-row sum cases, $p_{uv}$ is set to $\hat{q}$ (see also \S\S \ref{improvements}, \ref{other}).

In the Helmholtz equation, $c_{uv}$ is small for all $u,v$, indicating that all nodes are ``distant'' and that no single aggregation of the nodes will yield fast AMG convergence (indeed, it is known that multiple coarse grids are required \cite{wave_accuracy}).

\subsubsection{Aggregation Rules}
Ideal aggregates have strong internal connections and weaker external connections. To this end, we will be guided by four rules:
\ben
    \item Each node can be associated with one seed.
    \item A seed cannot be associated.
    \item Aggregate stronger affinitives before weaker.
    \item Favor aggregates with small energy ratios (cf. \S \ref{energy_correction}).
\een
Rules 1 and 2 prevent an associate from being transitively associated with multiple seeds. Otherwise, long chains of nodes might be aggregated together, creating aggregates with weak internal connections and very large energy ratios. Rule 3 favors strongly-connected aggregates. Rule 4 has dual purpose: (a) Ultimately, the energy ratio determines the AMG asymptotic convergence factor, hence this rule ensures good convergence. (b) Affinities are based on local information (relaxed TVs); their quantitative value becomes fuzzier as nodes grow apart \cite[\S 5]{alg_distance}. Since small energy ratios usually lead to small aggregates, affinities are indeed used only for {\it local} aggregation decisions. A typical coarsening ratio in our algorithm ranges between .3--.5.

\subsubsection{Aggregation Algorithm}
\label{agg_algorithm}
The main call $aggregate()$ gradually aggregates $\cN$ in $q$ stages, i.e., generates aggregate sets $S_1,\dots,S_q$ such that each $S_i$-aggregate is contained in some $S_{i+1}$-aggregate (Algorithm~\ref{aggregate}). The affinity threshold $\delta$ is monotonically decreased for stronger connections to be aggregated before weaker connections. Among the sets, we select that for which the coarsening ratio $\alpha:=|S_i|/n$ is closest to $\amax := g/\gamma$. This aims at bounding the total cycle work by $\approx 1 + \gamma \amax + (\gamma \amax)^2 + \dots \approx (1-g)^{-1}$ finest-level units, assuming the same fill-in at all levels (a more accurate $\alpha$ definition could be the ratio of coarse to fine edge numbers, which is easily tracked during aggregation).

$S$ is encoded by the $\status$ array: $\status(i)=\undecided:=-1$ denotes an undecided node, $\status(i)=\seed:=0$ a seed node, and $\status(i)>0$ indicates that $i$ is an associate of the seed $\status(i)$. High-degree nodes are apriori set to seeds, because a sun may only be another node's associate at a coarser level where its degree has dropped to $O(1)$. In each stage, we loop over undecided nodes and decide whether to aggregate each one with an existing seed, or with an undecided node that thereby becomes a new seed.

\begin{algorithm}
\caption{$S = aggregate(\bA, \bX, \amax, \Imax$)}
\label{aggregate}
\begin{algorithmic}[1]
    \STATE $\Imax \leftarrow 2$, $B(1..I_{\mathrm{max}}) \leftarrow \infty$, $n \leftarrow \sz(\bA)$.
    \STATE $n_c \leftarrow n$, $\alpha \leftarrow 1$, $i \leftarrow 0$, $\delta \leftarrow .9$.
    \STATE $\bC \leftarrow (c_{uv})_{u,v}$, where $c_{uv} \leftarrow$ (\ref{cuv}) with test vectors $\bX$, $\forall (u,v) \in \cE$.
    \STATE For each $u \in \cU$, $\status(u) \leftarrow \undecided, \aggsize(u) \leftarrow 1$.
    \STATE For each $u \in \{ u: |\cA_u| \geq 8 \mathrm{\;median}(\{\cA_v\}_v) \}, \status(u) \leftarrow 0$.
    \WHILE[Main aggregation loop]{($(\alpha \geq \amax$) and $(i < I_{\mathrm{max}})$)}
        \STATE $i \leftarrow i+1$, $\delta \leftarrow .6\delta$.
        \STATE $aggregationStage(\underline{status}, \underline{\bA}, \underline{n_c}, \underline{\bC}, \underline{\bX}, \underline{\aggsize}, \delta)$
        \STATE $\alpha \leftarrow n_c/n$, $B_i \leftarrow$ ($1-\alpha$ if $\alpha \leq \amax$, otherwise $1+\alpha$).
        \STATE[Save current aggregate set] $S_i \leftarrow \status$.
        \STATE For each $u \in \{v : S_i(v) = \seed \}, S_i(v) \leftarrow v$.
    \ENDWHILE
    \STATE $i \leftarrow \mathrm{argmin}(B)$.
    \STATE Return $S_i$.
\end{algorithmic}
\end{algorithm}

 $aggregationStage()$ relies on $bestSeed()$ to locate the closest seed of each node. Note: while $\bX$ and $\aggsize$ are modified in Algorithm~\ref{aggregationStage}, neither is subsequently used. This anticipates an improved $bestSeed()$ implementation in the next section that does utilize the updated values.
 
 Parameters in the code are passed by reference when underlined, otherwise passed by value.
 
\begin{algorithm}
\caption{$aggregationStage(\underline{status}, \underline{\bA}, \underline{n_c}, \underline{\bC}, \underline{\bX}, \underline{\aggsize}, \delta)$}
\label{aggregationStage}
\begin{algorithmic}[1]
    \STATE[Strong neighbors] $\bN \leftarrow \{ c_{uv} \geq \delta \max\left\{ \max_{s \not = u} c_{us}\,, \max_{s \not = v} c_{sv}\right\} \}$.
    \STATE $\bU \leftarrow \{u : \status(u) = \undecided \}$.
    \FOR[Undecided nodes with strong neighbors]{each $u$ in $\bU \cap \{u : \cA_u \cap \bN \not = \emptyset $}
        \IF[$i$'s status was changed earlier in this sweep]{$\status(i) \not = \undecided$}
            \STATE continue
        \ENDIF
        \STATE $s \leftarrow bestSeed(\bA, \bX, \aggsize, \bN, \bC, u)$.
        \IF[Identified seed neighbor, aggregate $u$ with $s$]{$s \not = \notfound$}
            \STATE $\status(s) \leftarrow \seed$, $\status(u) \leftarrow s$, $n_c \leftarrow n_c-1$.
            \STATE[Update TV values on new aggregate] For $k=1,\dots,K$, $x_{uk} \leftarrow x_{sk}$.
            \STATE $\aggsize(u),\aggsize(s) \leftarrow \aggsize(s)+1$.
        \ENDIF
    \ENDFOR
\end{algorithmic}
\end{algorithm}

\begin{algorithm}
\caption{$s \leftarrow bestSeed(\bA, \bX, \aggsize, \bN, \bC, u)$}
\label{bestSeed}
\begin{algorithmic}[1]
    \STATE[seed \& undecided $\delta$-affinitives] $S \leftarrow \bN \cap \{u : \status(u) \in \{\undecided, \seed\} \}$.
    \STATE {\bf if} $S = \emptyset$, return $\notfound$, {\bf else} [Closest neighbor] return $\argmax{v \in S} c_{uv}$. {\bf end if}
\end{algorithmic}
\end{algorithm}

\subsection{Energy-Corrected Coarsening}
\label{energy_correction}

\subsubsection{Energy Inflation}
\label{inflation}
\begin{figure}[htbp]
    \centering
    \begin{center}
    \begin{tabular}{cc}
      \includegraphics[height=1.5in]{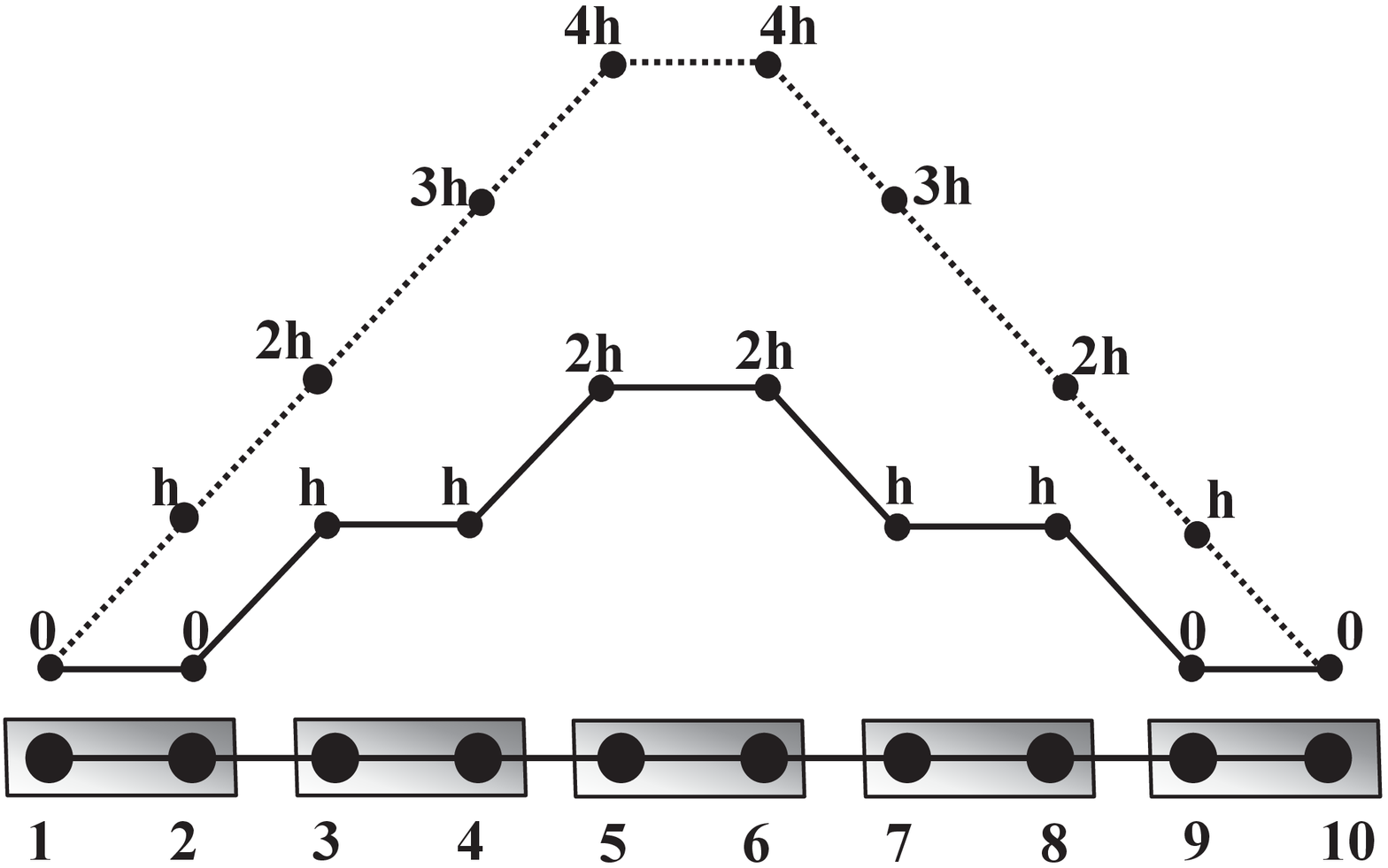}
      &
      \includegraphics[height=1.5in]{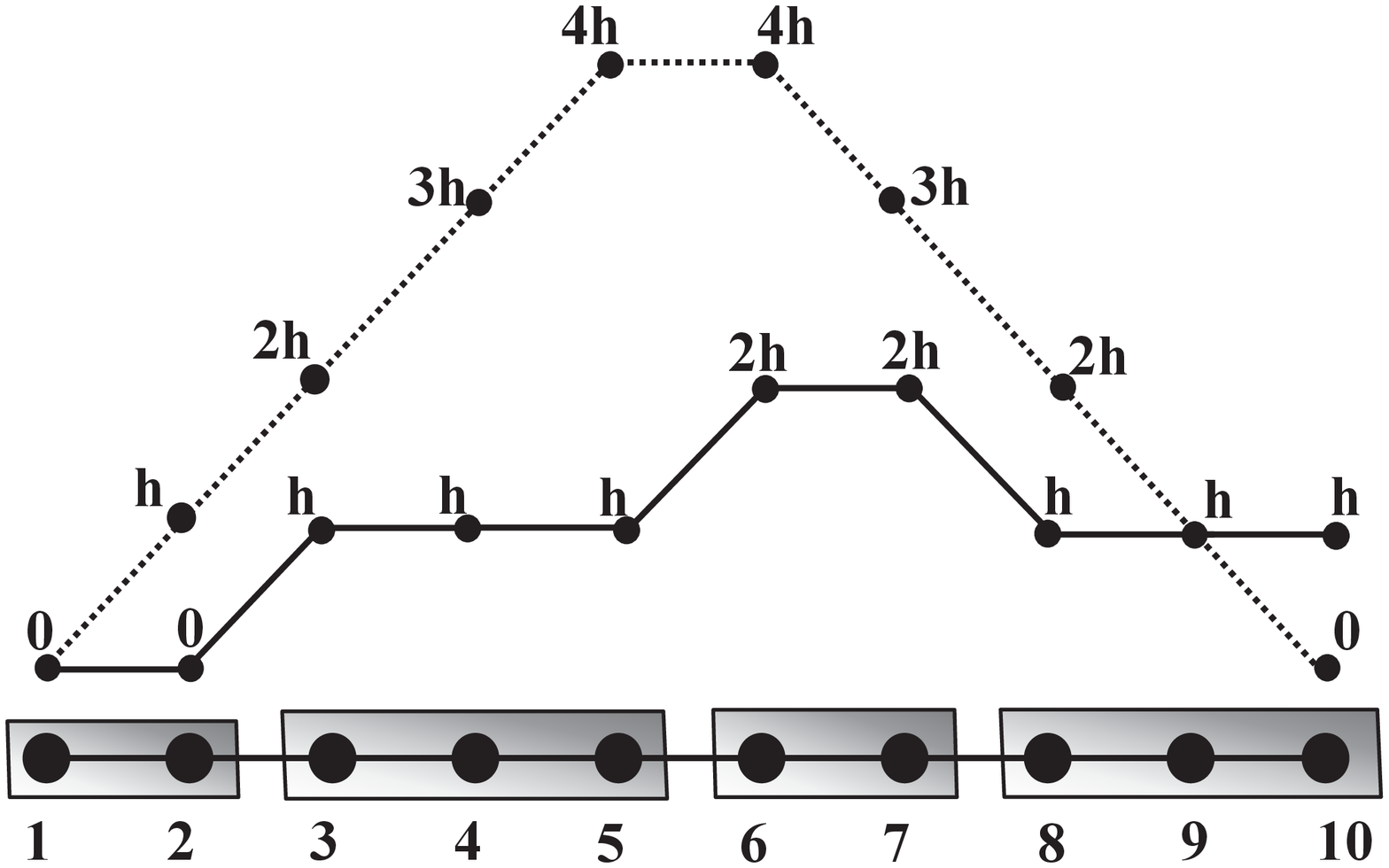}
      \\
      (a) & (b) \\
    \end{tabular}
    \end{center}
    \caption{The Galerkin correction $P \bee^c$ to a piecewise linear error $\bee$ in a path graph for (a) uniform 1:2 aggregation and (b) variable aggregate size. $h>0$ is arbitrary.}
    \label{path}
\end{figure}

The Galerkin coarse-level correction (\ref{min_ec}) $\bee^c$ is the best approximation to a smooth error $\bee$ in the energy norm. Braess \cite{braess} noted that this does not guarantee a good approximation in the $l_2$ norm. For example, if $\bee$ is a piecewise linear function in a path graph (1-D grid with $w \equiv 1$) coarsened by aggregates of size two, $P \bee^c$ is constant on each aggregate and matches $\bee$'s slope across aggregates, resulting in about half the fine-level magnitude; cf. Fig.~\ref{path}a.

An equivalent but more useful observation is that the energy of $\bP \bT \bee$ is twice larger than $\bee$'s, where $\bT \bee$ is a coarse representation of $\bee$, say,
\be
    \lp \bT e \rp_U := \frac{1}{|\cT_U|} \sum_{u \in \cT_U} \bee_u\,,\qquad U \in \cN^c\,.
    \label{t}
\ee
$\bT$ is called the {\it aggregate type} operator (alternatively, $\bT \bee$ could be defined as the vector of seed values. Either definition satisfies $\bT \bP = \bI^c$ and falls under the umbrella of a corresponding CR theory \cite[\S 2]{cr_james}). Rewrite (\ref{min_ec}) as
\be
   \mymin{\by^c} \left\{ \frac12 (\by^c)^{T} \bA^c \by^c - \by^c \bP^{T} \br \right\}\,,\qquad \br := \bA \bee\,.
   \label{min_ec_r}
\ee
For an ideal interpolation $\bP$ that satisfies $\bP \bT \bee = \bee$, (\ref{min_ec_r}) is minimized by $\by^c = \bT \bee$. A caliber-1 interpolation $\bP$ still satisfies $\bP \bT \bee \approx \bee$, but the first term in (\ref{min_ec_r}) is multiplied by the {\it energy inflation factor}
\be
    q(\bee) := \frac{E^c\left(\bT \bee\right)}{E\left(\bee\right)} = \frac{E\left(\bP \bT \bee\right)}{E\left(\bee\right)},
    \qquad E^c\left(\bee^c\right) := \frac12 \left(\bee^c\right)^{T} \bA^c \bee^c\,.
    \label{q}
\ee
Now (\ref{min_ec_r}) is minimized by $\by^c \approx q^{-1}\bT \bee$. As $\bee$ is not significantly changed by relaxation, its two-level ACF will be $\rho \approx 1 - 1/q$. In Fig.~\ref{path}a, $q \approx 2$ and $\rho \approx .5$.

\subsubsection{Energy Correction}
\label{mu}
Several energy inflation remedies may be pursued:
\bi
    \item[(A)] Increase the interpolation's caliber and accuracy. This leads to the fill-in troubles discussed in \S \ref{considerations}.
    \item[(B)] Accept an inferior two-level ACF of $1-1/q$ and increase the cycle index $\gamma$ to maintain it in a multilevel cycle. Unfortunately, not only does this increase complexity, the examples of \S \ref{grid_examples} demonstrate that $q$ can be arbitrarily large. The ACF cannot be improved by additional smoothing steps either, because it is governed by smooth mode convergence.
    \item[(C)] Correct the coarse level operator $\bA^c$ to match the fine level operator's energy during the setup phase.
    \item[(D)] Modify the coarse level correction $\bP \bee^c$ to match the fine level error $\bee$ during the solve phase.
\ei
In this section we consider option C. The Galerkin equation (\ref{galerkin}) is modified to
\be
    \bP^{T} \bA \bP \, \bee^c = \mu \bP^{T}\left(\bb - \bA \tbx\right)\,.
    \label{galerkin2}
\ee
The key question is how to choose $\mu$. (Note that if $\mu$ is constant, Options C and D are equivalent.) Motivated by Fig.~\ref{path}a and its two-dimensional analogue, Braess used $\mu=1.8$, but his V-cycle convergence for grid graphs was mesh-independent only if a fixed number of levels were used per cycle and if AMG was used as a PCG preconditioner. Moreover, no {\it predetermined global} factor exists that fits all error {\it corrections} in scenarios like Fig.~\ref{path}b, because the coarse-level solution depends on local inflation ratios, which vary among graph nodes.

On the other hand, a {\it local energy} correction $\mu$ does exist. Indeed, the quadratic energies are separable to nodal energies
\begin{subequations} \label{nodal}
\begin{eqnarray}
    E(\bx) &=& \sum_{u \in \cN} E_u(\bx)\,,\quad\quad\,\,
    E_u(\bx) := -\frac12 \sum_{v: v \not = u} a_{uv} \left(x_u-x_v\right)^2\,, \label{nodal_fine} \\
    E^c(\bx^c) &=& \sum_{U \in \cN^c} E^c_u(\bx^c)\,,\quad
    E^c_U(\bx^c) := -\frac12 \sum_{V: V \not = U} a^c_{UV} \left(x^c_U-x^c_V\right)^2\,, \label{nodal_coarse}
\end{eqnarray}
\end{subequations}
and define the {\it local inflation factor} as
\be
    q_U(\bx) := \frac{E^c_U\left(\bT \bx\right)}{\sum_{u \in U} E_u\left(\bx\right)}\,.
    \label{q_local}
\ee
In principle, a {\it local} $\mu_U$ can be designed using our TVs to at least partially offset $q_U$; unfortunately, new difficulties arise (cf. \S \ref{other_corrections}). Consequently, we chose to still scale all RHS entries by a flat $\mu$, but {\it modify the aggregation so that $q_U(\bx) \lessapprox Q$ for all smooth vectors $\bx$ and all $U \in \cN^c$}, where $Q>1$ is a parameter. Under this condition a global factor is effective, whose optimal value minimizes the overall convergence factor:
\be
    \mu = \argmin{\mu>0} \mymax{1 \leq q \leq Q} \left| 1 - \frac{\mu}{q} \right| =
        \argmin{\mu>0} \left\{ \left| 1 - \mu \right|, \left| 1 - \frac{\mu}{Q} \right| \right\} = \frac{2Q}{Q+1}\,.
    \label{optimal_mu}
\ee
In our algorithm we set $Q=2$ and $\mu=\frac43$. The expected smooth mode ACF is $Q/(Q+1)=\frac13$.

Alternatively, a {\it dynamic} (flat or local) $\mu$ can be computed for scaling the correction $P \bee^c$ during the solve phase; see \S \ref{adaptive}.

\subsubsection{Grid Examples}
\label{grid_examples}
\begin{figure}[htbp]
    \centering
    \begin{center}
    \begin{tabular}{cc}
      \includegraphics[height=1.17in]{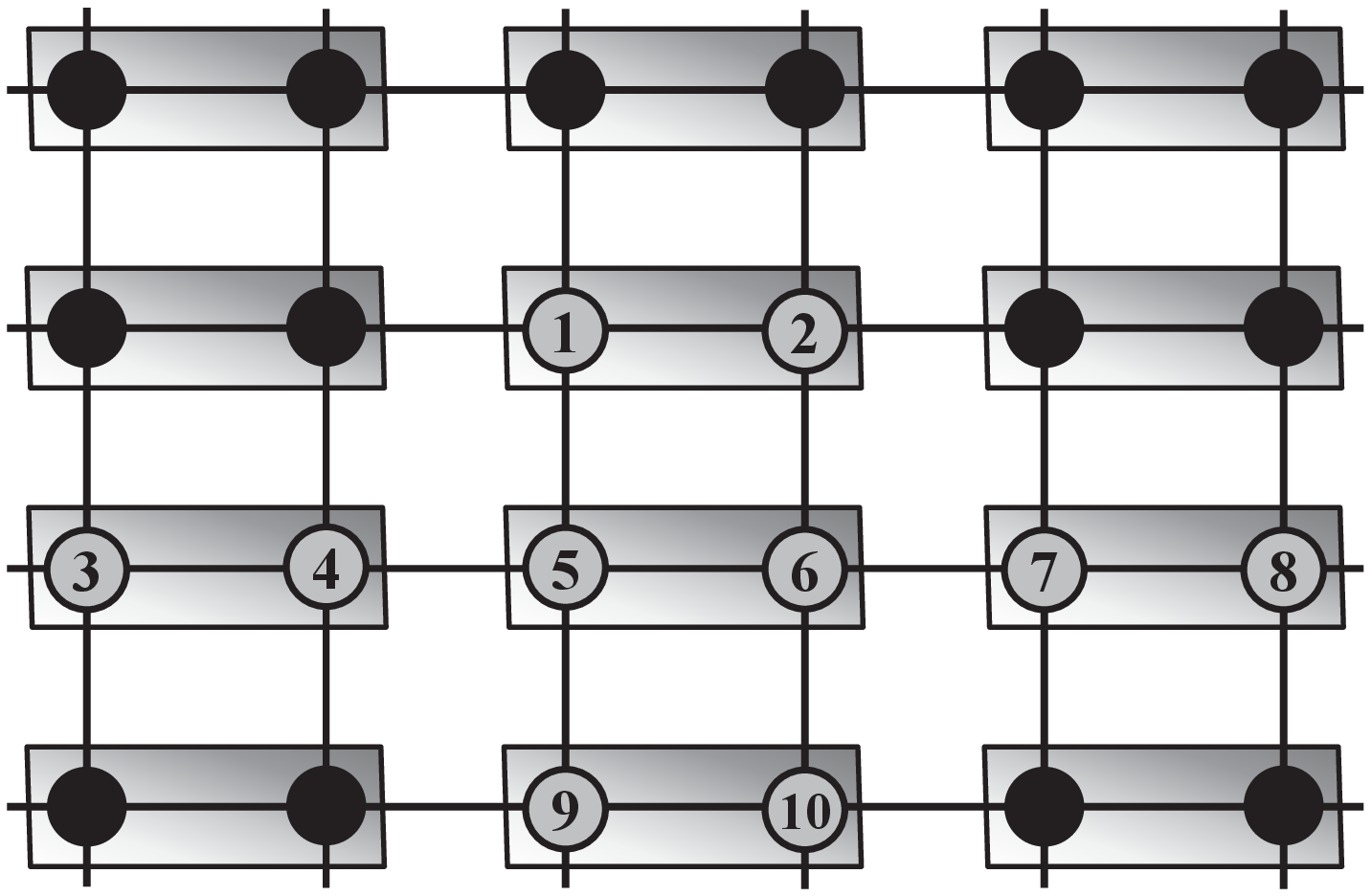}
      &
      \includegraphics[height=1.17in]{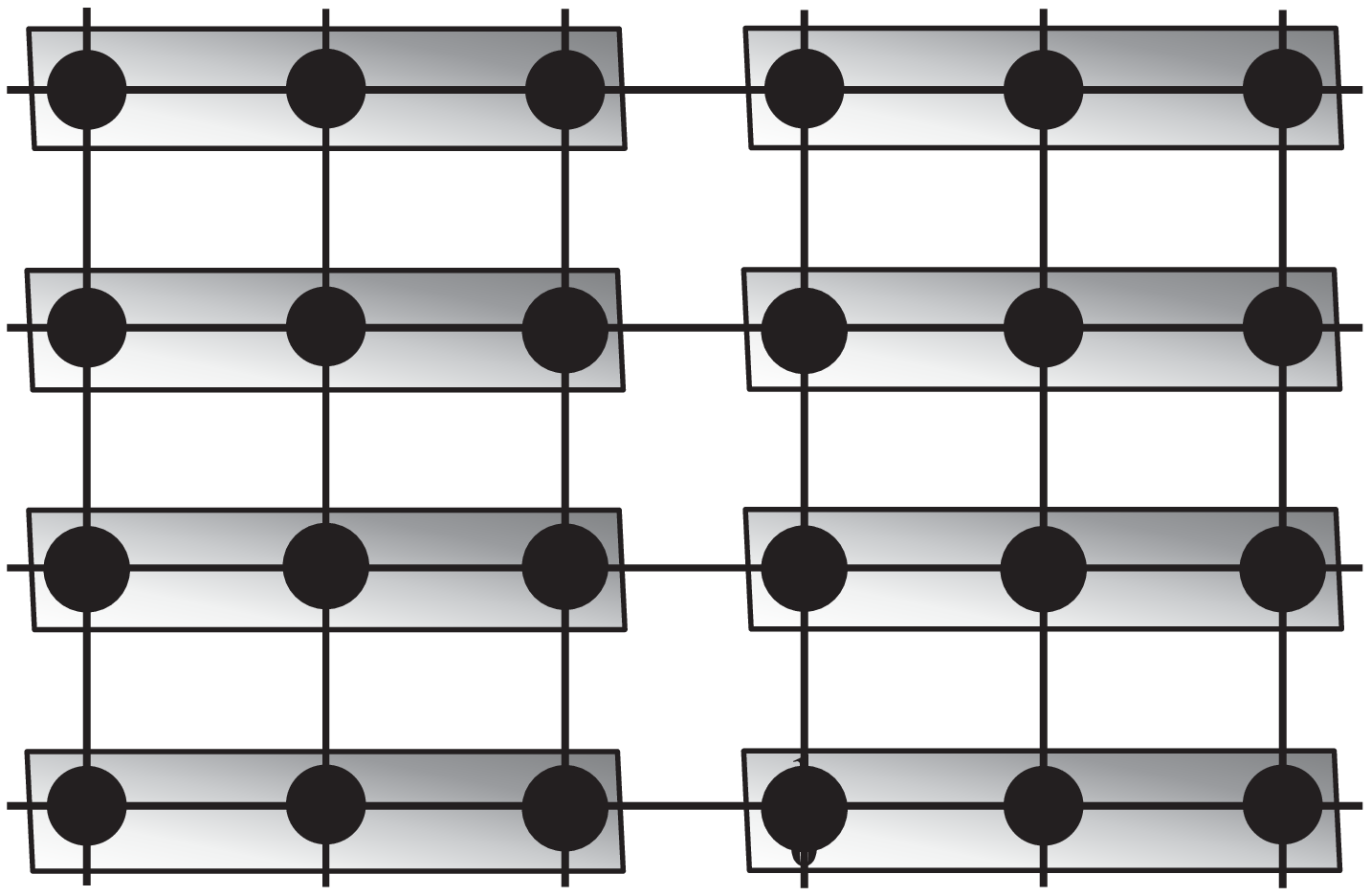}
      \\
      (a) & (b) \\
      \\
      \includegraphics[height=1.17in]{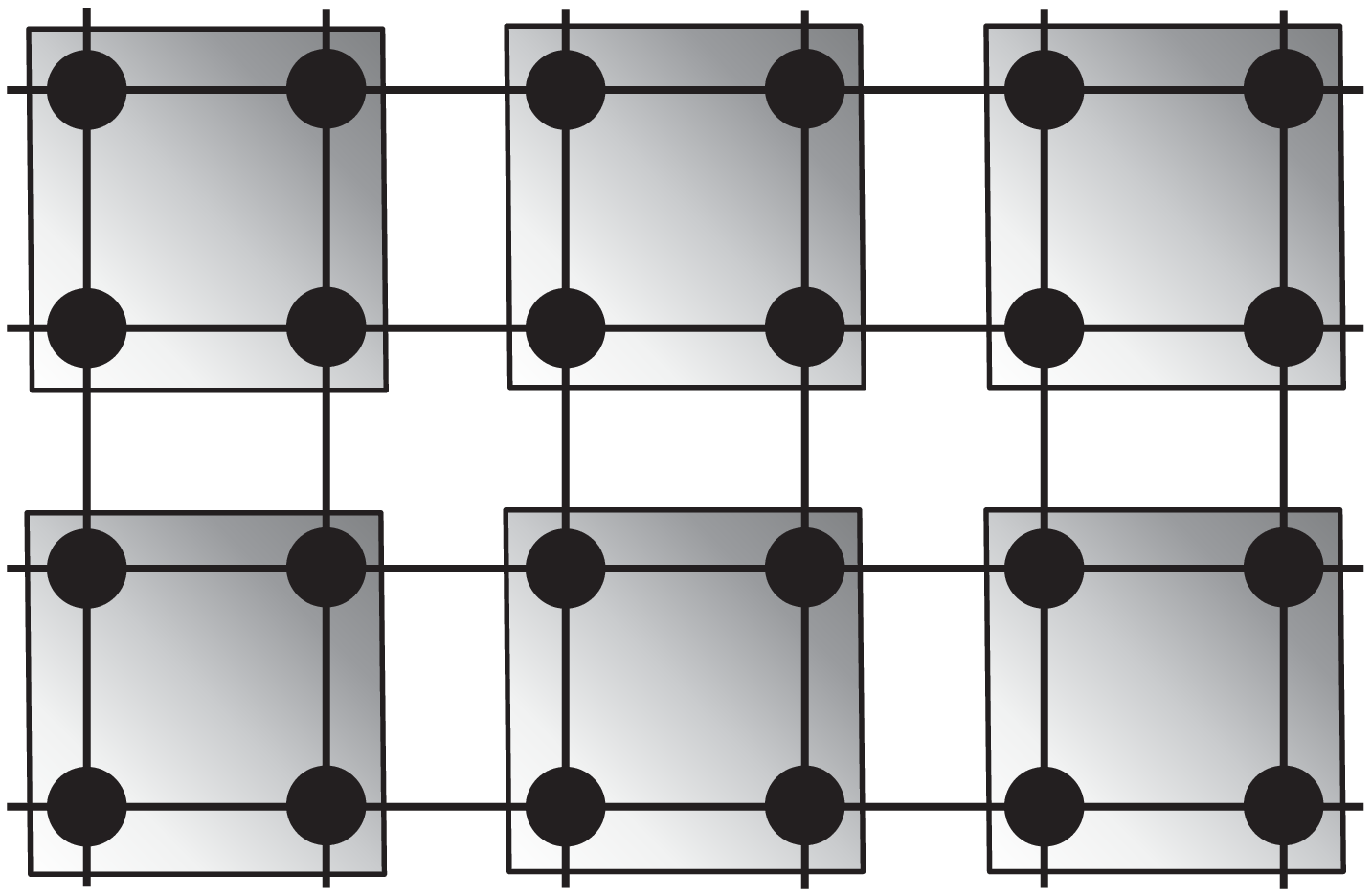}
      &
      \includegraphics[height=1.17in]{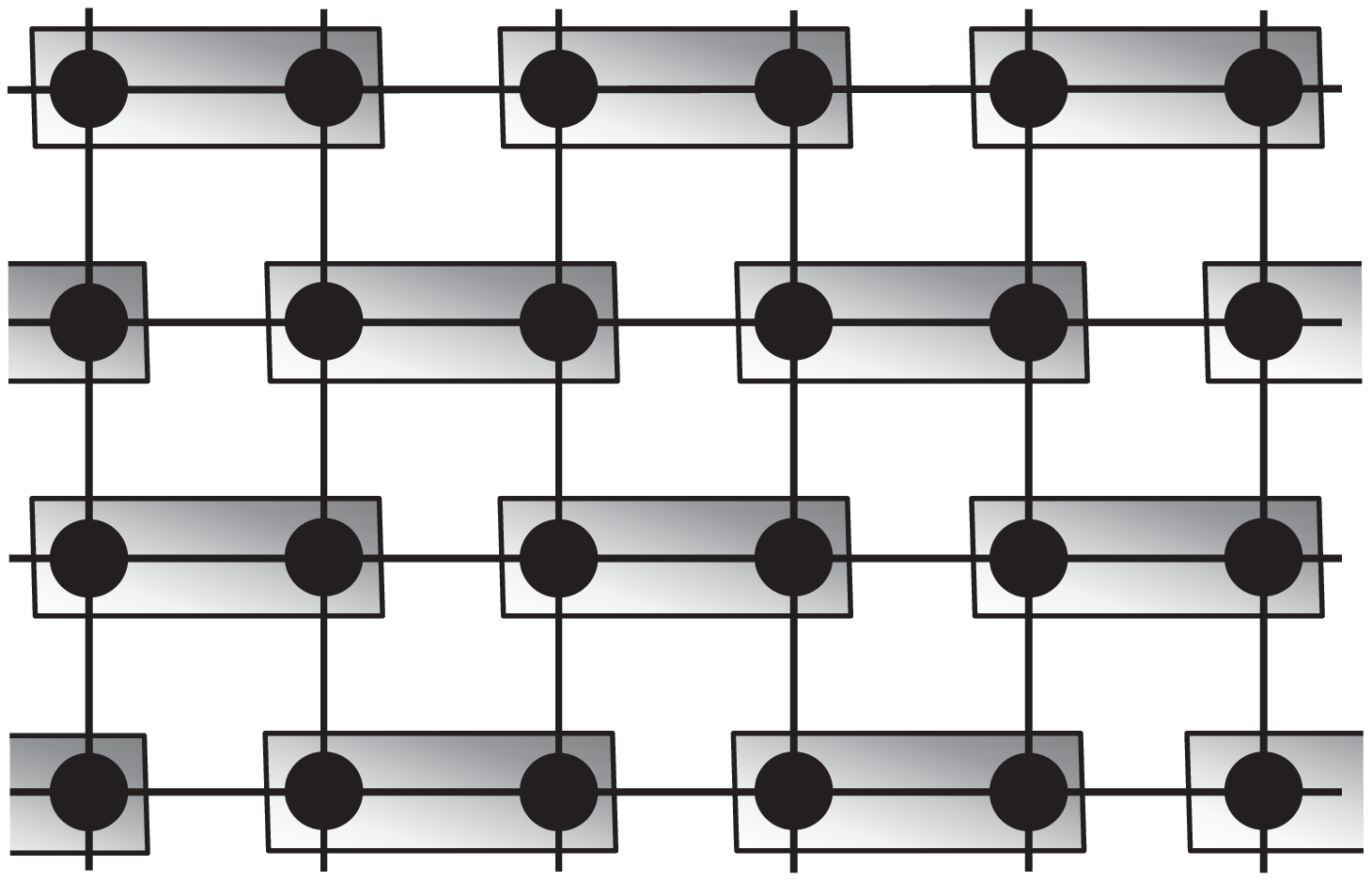}
      \\
      (c) & (d) \\
    \end{tabular}
    \end{center}
    \caption{2-D grid coarsening patterns: (a) 1:2 semi-coarsening ($Q=2$). (b) 1:3 semi-coarsening ($Q=3$).
    (c) Full coarsening ($Q=2$). (d) Staggered semi-coarsening ($Q=3$).}
    \label{grid2d}
\end{figure}

The local energy inflation varies considerably with aggregate size, shape and alignment. Consider the unweighted 2-D grid graph whose nodes are the 2-D coordinates $\{(i_1,i_2): i_1,i_2 \in \mathbb{Z}\}$ in the plane $t_1$-$t_2$. Since all affinities are equal, the aggregation algorithm of may locally create archetype constellations such as Fig.~\ref{grid2d}a--d depending on node ordering during the aggregation stages.

Low-energy vectors $\bx$ are locally linear in $t_1$ and $t_2$, so their nodal values are modelled by $\bx(t_1,t_2)=a_1 t_1+a_2 t_2$. In Fig.~\ref{grid2d}a, the worst local energy ratio $Q$ at any aggregate, say $U = \{5,6\}$, is
\small
\begin{eqnarray*}
    E_5(a_1,a_2) &=& \frac12 \sum_{v \in \{1,4,6,9\}} (\bx_5-\bx_v)^2 = 2 a_1^2 + 2 a_2^2\,, \\
    E_6(a_1,a_2) &=& \frac12 \sum_{v \in \{2,5,7,10\}} (\bx_5-\bx_v)^2 = 2 a_1^2 + 2 a_2^2\,, \\
    E^c_U(a_1,a_2) &=& \frac12 \left\{ \left(.5(\bx_5+\bx_6)-.5(\bx_3+\bx_4)\right)^2 +
                   \left(.5(\bx_5+\bx_6)-.5(\bx_7+\bx_8)\right)^2 \right. + \\
                   && \left. 2 \left(.5(\bx_5+\bx_6)-.5(\bx_1+\bx_2)\right)^2 +
                   2 \left(.5(\bx_5+\bx_6)-.5(\bx_9+\bx_{10})\right)^2 \right\} \\
                   &=& 4 a_1^2 + 2 a_2^2\,, \\
    q_U(a_1,a_2) &=& \left(4 a_1^2 + 2 a_2^2\right)/\left(2 a_1^2 + 2 a_2^2\right) \quad \Longrightarrow \quad
    Q = \mymax{a_1,a_2} q_U(a_1,a_2) = q_U(1,0) = 2\,.
\end{eqnarray*}
\normalsize
Similar analysis yields $Q=3$ in Fig.~\ref{grid2d}b and $Q=2$ in Fig.~\ref{grid2d}c. Small energy ratios thus occur when aggregates are small (consisting of 2--3 nodes each) and not elongated (a size-4 aggregate whose longest ``side'' is $2$ is preferable to chain of size $3$).

We might conjecture that aggregates of size $2$ as in Fig.~\ref{grid2d}a always ensure $Q \approx 2$, but this turns out to be a fallacy: the staggered coarsening in Fig.~\ref{grid2d}d manifests energy ratios as high as $Q=3$. Even worse, its $d$-dimensional analogue yields $Q=d+1$ -- unbounded as $d$ increases! Such aggregate constellations must be avoided. Fortuitously, we already possess the tool to signal and repel them: test vectors.

\subsubsection{Modified Aggregation Algorithm}
\label{agg_algorithm2}

The decision on aggregating $u$ with a seed $s$ during an aggregation stage (Algorithm~\ref{aggregationStage}) will now be based on both affinities and local TV energy ratios. $s$ is still required to be a $\delta$-affinitive neighbor of $u$. Additionally, for each TV we compare the nodal energy $E_u$ before and after aggregation, and aggregate only if the ratio is sufficiently small for all TVs.

Note that the nodal energy (\ref{nodal_fine}) is a quadratic in $x_u$ and $\{x_v\}_{v \in \cA_u}$. We define the more general functional
\be
    E_u(\bx;y) := \frac12 a_{uu} y^2 - B_u(\bx) y + C_u(\bx)\,, B_u(\bx) := \sum_{v \in \cA_u} w_{uv} x_v\,, C_u(\bx) := \frac12 \sum_{v \in \cA_u} w_{uv} x_v^2\,,
    \label{nodal_quad}
\ee
so $E_u(\bx) = E_u(\bx;x_u)$. The energy inflation is estimated by
\be
    q_{ut} := \mymax{1\leq k \leq K} \frac{E\left(\bx^{(k)};x_t\right)}{E\left(\bx^{(k)};\overline{x}_u\right)}\,,\,
    \overline{x}_u := \argmin{y} E_u(\bx;y) = \frac{B_u}{a_{uu}} = \frac{\sum_{v \in \cA_u} w_{uv} x_v}{\sum_{v \in \cA_u} w_{uv}}\,.
    \label{qut}
\ee
The numerator of $q_{ut}$ is the post-aggregation energy. The denominator is the smallest possible local energy of the fine-level TV. This aims at reducing the high-energy modes in $x^{(k)}$ that may obliterate the inflation estimate, and is equivalent to a temporary Jacobi relaxation step at $u$, which is also popular in recent bootstrap AMG works \cite[\S 2.2]{bamg}. (A block-relaxation step at $u$ and $s$ could instead be applied to increase accuracy at a higher cost, but we have not pursued it in the ``lean'' spirit of LAMG.)

$bestSeed()$ is replaced with a new implementation (Algorithm~\ref{bestSeed2}): $\{q_{ut}\}_t$ are computed for all $\delta$-affinitive seed neighbors $t$ of $u$. Among the neighbors for which $q_{ut} \leq 2.5$, $u$ is aggregated with that which represents the smallest aggregate (if the minimum is attained by multiple $t$'s, one of them is arbitrarily chosen). If all $q_{ut} > 2.5$, $u$ is not aggregated at all. Slightly larger ratios than the target $Q=2$ are accepted because TVs also contain small-magnitude high-energy modes for which strict ratios are neither necessary nor attainable. This especially applies at coarse levels, where more TVs are used and the chance of increasing the maximum TV energy ratio rises.

\begin{algorithm}[h]
\caption{$s \leftarrow bestSeed(\bA, \bX, \aggsize, \bN, \bC, u)$}
\label{bestSeed2}
\begin{algorithmic}[1]
    \STATE[seed \& undecided $\delta$-affinitives] $S \leftarrow \bN \cap \{u : \status(u) \in \{\undecided, \seed\} \}$.
    \STATE {\bf if} {$S = \emptyset$}, $return \notfound$; {\bf end if}
    \STATE For each $k=1,\dots,K$, compute $A_u(\bx^{(k)}), B_u(\bx^{(k)})$ using (\ref{nodal_quad}).
    \STATE For each $k=1,\dots,K$, $F_k \leftarrow E_u(\bx^{(k)};B_u(\bx^{(k)})/a_{uu})$ using (\ref{horner}).
    \STATE For each $k=1,\dots,K$, $t \in S$, $C^{(k)}_t \leftarrow E_u(\bx^{(k)};x^{(k)}_t)$ using (\ref{horner}), $q_t \leftarrow C^{(k)}_t/F_k$.
    \STATE[Neighbors with small energy inflation] $\tilde{S} \leftarrow \{ t \in S : q_t \geq 2.5 \}$.
    \STATE {\bf if} ($S = \emptyset$), return $\notfound$, {\bf else} return $\argmin{v \in \tilde{S}} \aggsize(v)$. {\bf end if}
\end{algorithmic}
\end{algorithm}

To save work, for each TV $\bx^{(k)}$ we first compute the terms $B_u$ and $C_u$ in (\ref{nodal_quad}), and subsequently evaluate the quadratic $E_u(\bx^{(k)};y)$ for each $y = x^{(k)}_v, v \in \cA_u$ and $y = \overline{x}^{(k)}_u$. Horner's rule \cite[p.~8]{poly_book} is applied to save a multiplication:
\be
    E_u(\bx^{(k)};y) = \left(\frac12 a_{ii} y - B_u\left(\bx^{(k)}\right) \right) y + C_u\left(\bx^{(k)}\right)\,.
    \label{horner}
\ee
The complexity of $bestSeed()$ is therefore at most $O(K |\cA_u|)$.

Two advantages of the low-degree elimination (\S \ref{elimination}) are (a) largely preventing worst-case 1-D scenarios such as Fig.~\ref{path}b, where it is impossible to obtain low energy ratios without excessively increasing the coarsening ratio; and (b) increasing the number of neighbors of $u$ and the chance of locating a seed $t$ with small energy inflation.

\subsubsection{Connected Component Assembly}
\label{connected_comp}
Linear-time algorithms for identifying the connected components of a graph exist \cite{tarjan, sedgewick}, and could be applied prior to LAMG invocation to reduce $G$ to a singly-connected graph. This is optional, as LAMG naturally represents components as the disconnected nodes of each level.

\begin{figure}[htbp]
    \centering
    \begin{center}
      \includegraphics[height=1.8in]{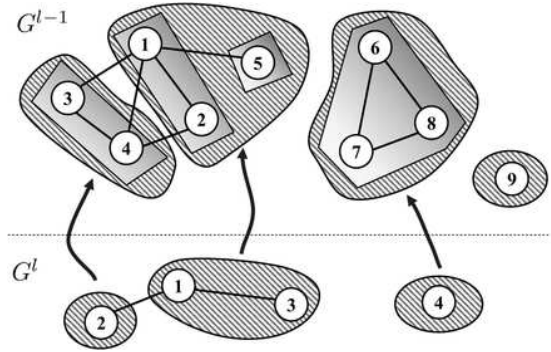}
    \end{center}
    \caption{The three components of $G^l$ (\{2\}, \{1,3\} and \{4\}) are interpolated to the $G^{l-1}$ components \{3,4\}, \{1,2,5\} and \{6,7,8\}. $G^{l-1}$'s fourth component is \{9\}.
    }
    \label{component_interp}
\end{figure}

We start by computing the coarsest graph $G^L$'s components $\{\cN^L_m\}_{m=1}^{M^L}$, e.g., using Tarjan's algorithm \cite{tarjan}. Let $\bP^L$ be the $\kth{L}$ level interpolation and $\{\cT^L_U\}_{U \in \cN^L}$ be the corresponding aggregate set (cf. \S \ref{aggregation}). The {\it set interpolation} $\bP^L \cU$ of a set $\cU^L \subseteq \cN^L$ is the support of its characterstic function's interpolant,
\be
    \bP^L \cU^L := \bigcup_{U \in \cU^L} \cT^L_U\,.
    \label{set_interp}
\ee
An $A^L$-edge between $U$ and $V$ cannot exist if $\cT^L_U$ and $\cT^L_V$ are disconnected. Thus, $G^{l-1}$'s connected components are the interpolants $\bP^L \cN^L_m$, $m=1,\dots,M^L$, plus any disconnected $G^{l-1}$ nodes. All $l-1$ components are interpolated to level $l-2$ and amended with disconnected $G^{l-2}$ nodes, and so on. When $l=1$ is reached, we obtain the set of connected components of $G^1=G$. See Fig.~\ref{component_interp} and Algorithm~\ref{connectedComponents}.

\begin{algorithm}
\caption{$N \leftarrow connectedComponents(\{G^l, \bP^l\}_{l=1}^L)$}
\label{connectedComponents}
\begin{algorithmic}[1]
    \STATE $N \leftarrow $ connected components of $G^L$, using Tarjan's algorithm \cite{tarjan}.
    \FOR{$l = L,L-1,\dots,2$}
        \STATE For each $\cU \in N$, $\cU \leftarrow \bP^l \cU$.
        \IF{$l$ is an \elim level}
            \FOR{each elimination stage (cf. \S \ref{elimination}) $i=q,q-1,\dots,1$}
                \STATE $\cZ_i \leftarrow$ be the set of disconnected nodes this stage.
                \STATE $N \leftarrow N \cup \left\{\bP^l_{i-1} \cdot \ldots \cdot \bP^l_1 \cZ_i \right\}$.
            \ENDFOR
        \ENDIF
    \ENDFOR
\end{algorithmic}
\end{algorithm}

\subsection{Solve Phase}
\label{solve}

The solve phase consists of multigrid cycles \cite[\S~1.4]{guide}. Each $l<L$ is assigned a cycle index $\gamma^l$ and pre- and post-relaxation sweep numbers $\nu^l_1,\nu^l_2$. If $l+1$ is an \elim level, $\gamma^l=1$ and $\nu^l_1=\nu^l_2=0$; otherwise,
\be
    \gamma^l :=
    \begin{cases}
        \gamma\,,& |\cE^l| > .1 |\cE|\,, \\
        \min\left\{2, g |\cE^{l+1}|/|\cE^l| \right\}\,,& \text{otherwise} \,,
    \end{cases}\qquad
    \nu_1^l = 1\,,\quad\nu_2^l = 2\,.
    \label{cycle_index}
\ee
At fine levels, $\gamma=1.5$ is employed; this value is theoretically marginal for a bounded multilevel ACF if the smoothest errors' two-level convergence factor is $\approx \frac13$ \cite[\S 6.2]{guide}, which seems to be implied by (\ref{optimal_mu}). Notwithstanding, worst-case energy ratios occur infrequently in practice. This issue is further diminished by the adaptive energy correction of \S \ref{adaptive}. At coarse levels, $\gamma^l$ is increased to maximize error reduction while incurring a bounded work increase. The total cycle work is about $3/(1-g)$ relaxation sweeps (cf. \S \ref{agg_algorithm}), so $g=.7$ seems like a reasonable value.

Three GS sweeps per level provide adequate smoothing, especially in light of the coarse-level correction's crudeness. It may be possible to escape with fewer sweeps, but the relaxation work is anyway dominated by the setup time, as illustrated in \S \ref{smorgasbord}.

Note that we use fixed settings for all cycle parameters: no optimization or fine tuning is required for specific graphs.

If level $2$ is an \elimnospace, $\bb$ is restricted to $\bb^2$ only once, and all cycles are applied to $\bA^2 \bx^2 = \bb^2$, followed by interpolating the final $\bx^2$ to $\bx$.

\subsubsection{Adaptive Energy Correction}
\label{adaptive}

Instead of fixing $\mu$ (cf. \S \ref{energy_correction}), one can modify the correction to smooth errors during the cycle. Let $l$ be any level such that $l+1$ is an \agg level. Whenever we transition to level $l$ from level $l-1$, $\vartheta$ sub-cycles are applied to $\bA^l \bx^l = \bb^l$, where $\vartheta$ is $1$ or $2$. We save the iterants $\bx^l_i$ obtained after the pre-relaxation of each sub-cycle. Before switching back to level $l-1$, the final iterant $\bx^l$ is replaced with $\by^l$, where
\be
    \by^l = \bx^l + \balpha_1 \left( \bx^l_1 - \bx^l \right) + \cdots + \balpha_{\vartheta} \left( \bx^l_{\vartheta} - \bx^l \right)
    \label{recomb}
\ee
and $\{\alpha_i\}_{i=1}^{\vartheta}$ are chosen such that $\|\bb^l - \bA^l \by^l\|_2$ is minimized (this is an $n_l \times \vartheta$ least-squares problem that is solved in $O(n_l)$ time). This iterant recombination \cite[\S 7.8.2]{trot} diminishes smooth errors $\bx^l_i-\bx^l$ that were not eliminated by $\kth{(l+1)}$-level corrections. Since the initial {\it residuals} obtained after interpolation from level $l+1$ are {\it not} smooth, a residual minimization is only effective after $\bx^l_i - \bx^l$ is smoothed. To maximize iterant smoothness, more relaxation sweeps are performed {\it after} returning from a coarse level, hence $\nu^l_1=1, \nu^l_2=2$ in (\ref{cycle_index}) (some pre-smoothing is still needed: this choice was superior to $\nu^l_1=0, \nu^l_2=3$ in numerical experiments).

This acceleration is superior to CG because it is performed at multiple levels. Iterant recombination at coarse levels was first suggested by A. Brandt and has been recognized as an effective multigrid tool \cite[Remark~7.8.5]{trot}. In LAMG, recombination occurs more frequently at coarser levels (since $\gamma>1$), where it is less expensive.

\begin{figure}[htbp]
    \centering
    \begin{center}
      \includegraphics[height=1.85in]{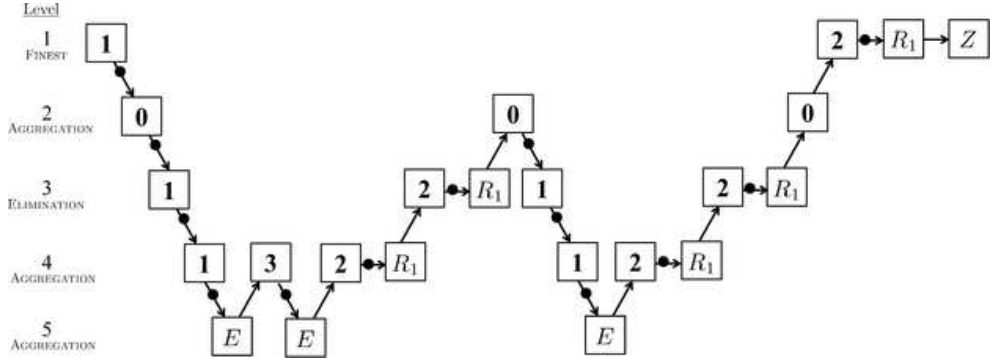}
    \end{center}
    \caption{A four-level cycle. A squared number denotes a number of GS sweeps; $E$ denotes an exact solver (\S \ref{coarsest}); $Z$ denotes the Grahm-Schmidt loop (\ref{grahm}); down-arrows denote RHS restrictions (\ref{elim_coarse}) or (\ref{galerkin2}); up-arrows denote corrections (\ref{xz})--(\ref{xf}) or (\ref{correction}). Iterants are saved at the black dots before the restriction; $R_{\vartheta}$ denotes a $(\vartheta+1)$-iterant recombination (\ref{recomb}).}
    \label{cycle}
\end{figure}

\subsubsection{Zero-Mode Orthogonalization}
\label{orthog}
The null-space components of $\bx$ are not determined by the cycle. At the end of each cycle, we add the Grahm-Schmidt procedure
\be
    \text{For } m=1,\dots,M\,,\qquad \bx \leftarrow \bx + \frac{\alpha_m - \bu_m^{T} \bx}{\bu_m^{T} \bu_m} \, \bu_m\,.
    \label{grahm}
\ee
Note that the zero modes $\bu_1,\dots,\bu_M$ are known once the connected components are assembled. (Alternatively, the cycle could be reformulated using the full approximation scheme. The constraints (\ref{xcompat}) would be transferred to the coarsest level and solved concurrently with (\ref{system}) \cite[\S 5.6]{guide}.)

\subsubsection{Coarsest Solver}
\label{coarsest}

If relaxation converges fast at a certain level, the linear system is iteratively solved to sufficient accuracy by GS solve iterations (\S \ref{relax}). If it is slow, the coarsest graph $G^L$ should be small enough for a direct solver to be fast. This size is implementation-specific; in our program, coarsening terminates when $n^L \leq 150$ and the augmented system
\be
    \begin{pmatrix}
    \bA^L & \bU^L \\
    \left(\bU^L\right)^{T} & \bzero
    \end{pmatrix}
    \begin{pmatrix}
    \bx^L \\
    \bxi^L
    \end{pmatrix}
    =
    \begin{pmatrix}
    \bb^L \\
    \bzero
    \end{pmatrix}
    \label{coarsest_system}
\ee
is directly solved. $\bU^L$ is a matrix whose columns are the characteristic functions of the connected components $\cN^L_1,\dots,\cN^L_{M^L}$ of $G^L$. The dual vector $\bxi^L$ is discarded.

\section{Numerical Results}
\label{results}
Next, we provide supporting evidence for LAMG's practical efficiency for a wide range of graphs.

\subsection{Smorgasbord}
\label{smorgasbord}

An unoptimized LAMG object oriented \matlab 7.10.0 (R2010a) serial implementation was developed and tested on a diverse set of \numgraphs real-world graphs with up to \maxedges edges, collected from the following sources:
\bi
    \item The University of Florida Sparse Matrix collection (UF) \cite{uf_collection}.
    \item C. Walshaw's graph partitioning archive \cite{walshaw}.
    \item I. Safro's MLogA results archive at Argonne National Laboratory \cite{mloga}.
    \item The FTP site of the DIMACS Implementation Challenges \cite{dimacs}.
\ei
The graphs originated from a plethora of applications: CFD airplane and car FEM meshes; RF electrical circuits; combinatorial optimization; model reduction benchmarks; social networks; web page networks; and many others. Experiments were performed on a 64-bit Windows 7 Dell Inspiron 580 (3.2 GHz CPU; 8GB RAM).

For each graph, a compatible RHS was generated by identifying a pair of nodes $s$, $t$ in the same connected component and setting $b_s=1$, $b_t=-1$ and $b_u=0$, $u \not \in \left\{s,t\right\}$ ($\bA \bx = \bb$ models an effective capacitance problem \cite[\S 2]{icm10}). LAMG setup was executed with two aggregation stages at each level. That is, $\delta=.9$ was set in the first stage of Algorithm~\ref{aggregate}, and if a second stage was performed, it used $\delta=.54$. The solve phase was then invoked twice, with a flat $\mu=\frac43$ energy correction and an adaptive correction; each solve started from a initial random guess and proceeded till the residual $l_2$-norm was reduced by $10^8$. Five performance measures were computed for each graph:
\bi
    \item {\it Setup time per edge $\tsetup$} [seconds].
    \item {\it Solve time per edge per significant figure $\tsolve$} [seconds], for the adaptive scheme. If residual norm after $i$ iterations was $r_i$ and $p$ iterations were executed, $\tsolve := t/(m \log_{10}(r_0/r_p))$, where $t$ was the total solve time.
    \item {\it Total time per edge $\ttotal = \tsetup + 10 \tsolve$} [seconds].
    \item {\it Asymptotic Convergence Factor (ACF)} of the flat and adaptive schemes, estimated as $(r_p/r_0)^{1/p}$.
    \item {\it Adaptive correction gain $A$}: the ratio of flat-to-adaptive solve times.
\ei

\begin{figure}[htbp]
    \centering
    \begin{center}
    \begin{tabular}{cc}
      \includegraphics[height=1.85in]{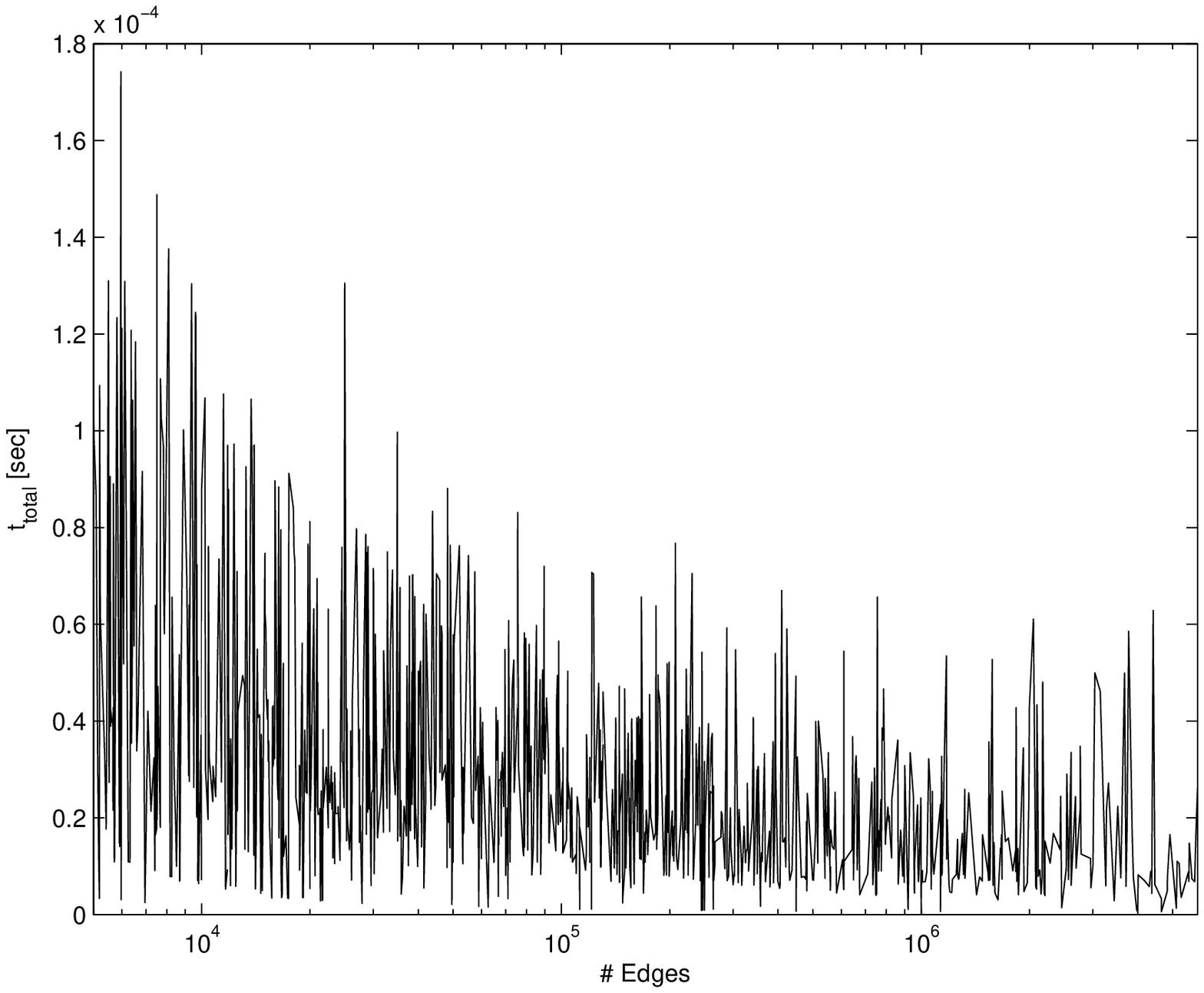}
      &
      \includegraphics[height=1.85in]{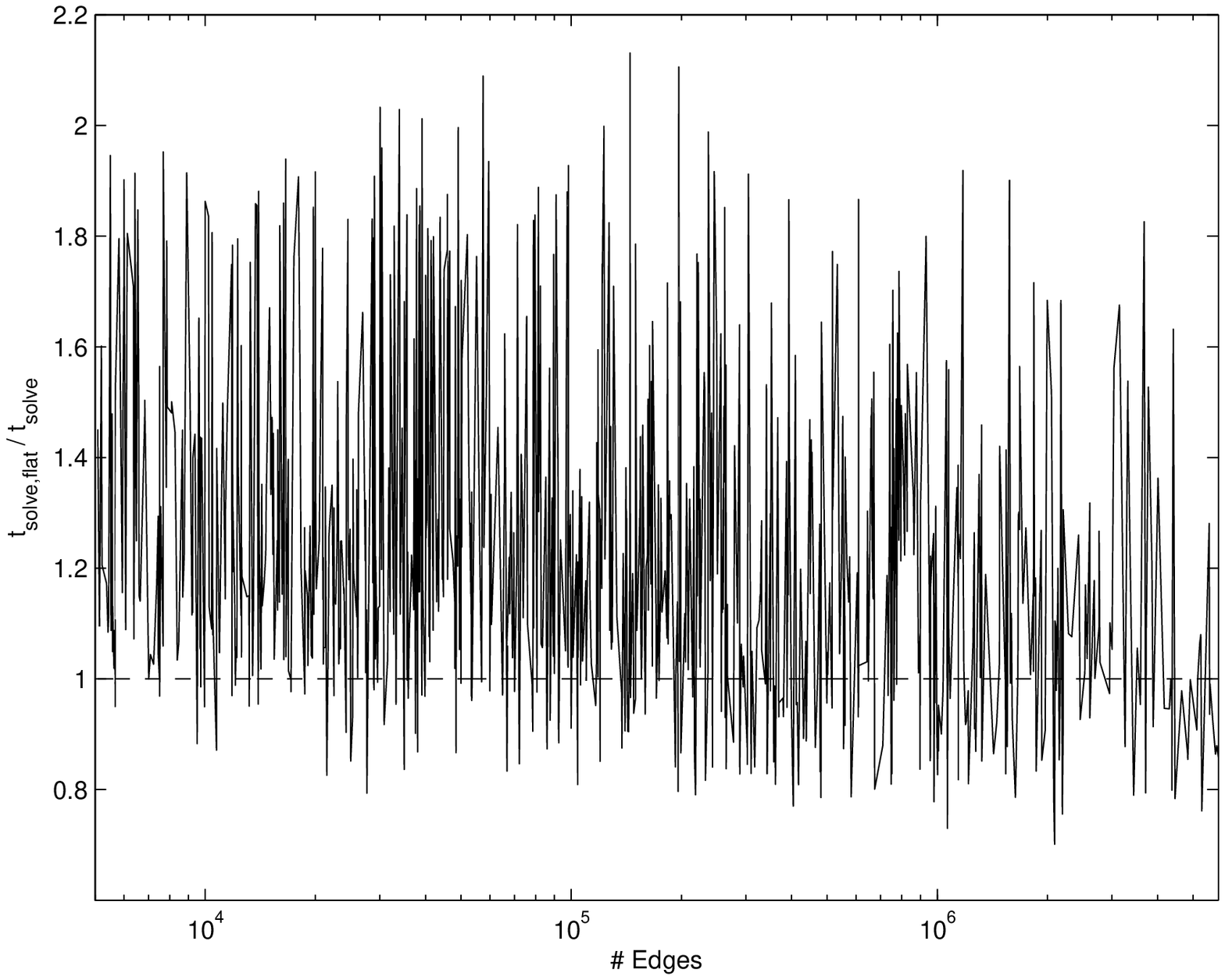}
      \\
      (a) & (b) \\
      \includegraphics[height=1.85in]{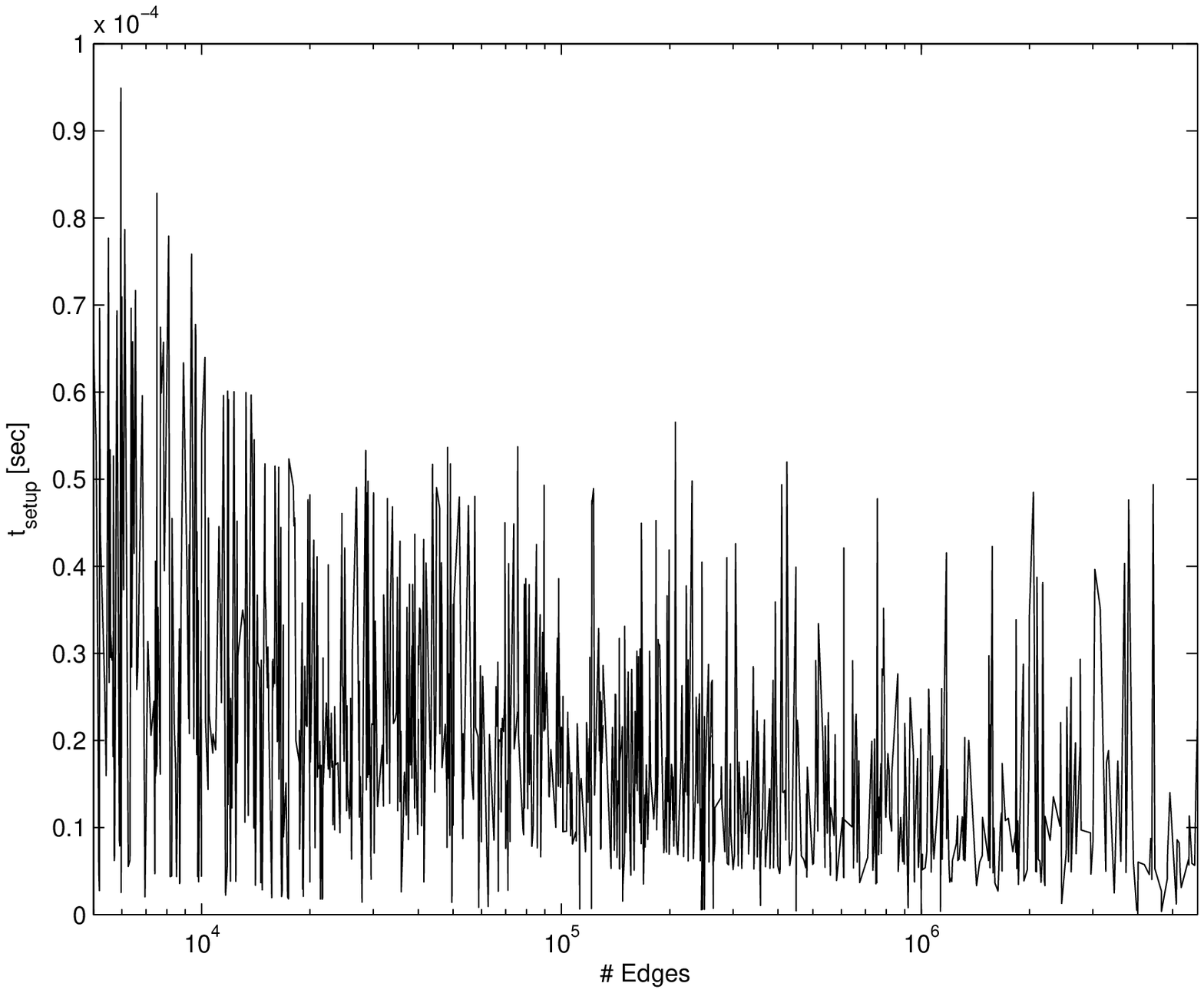}
      &
      \includegraphics[height=1.85in]{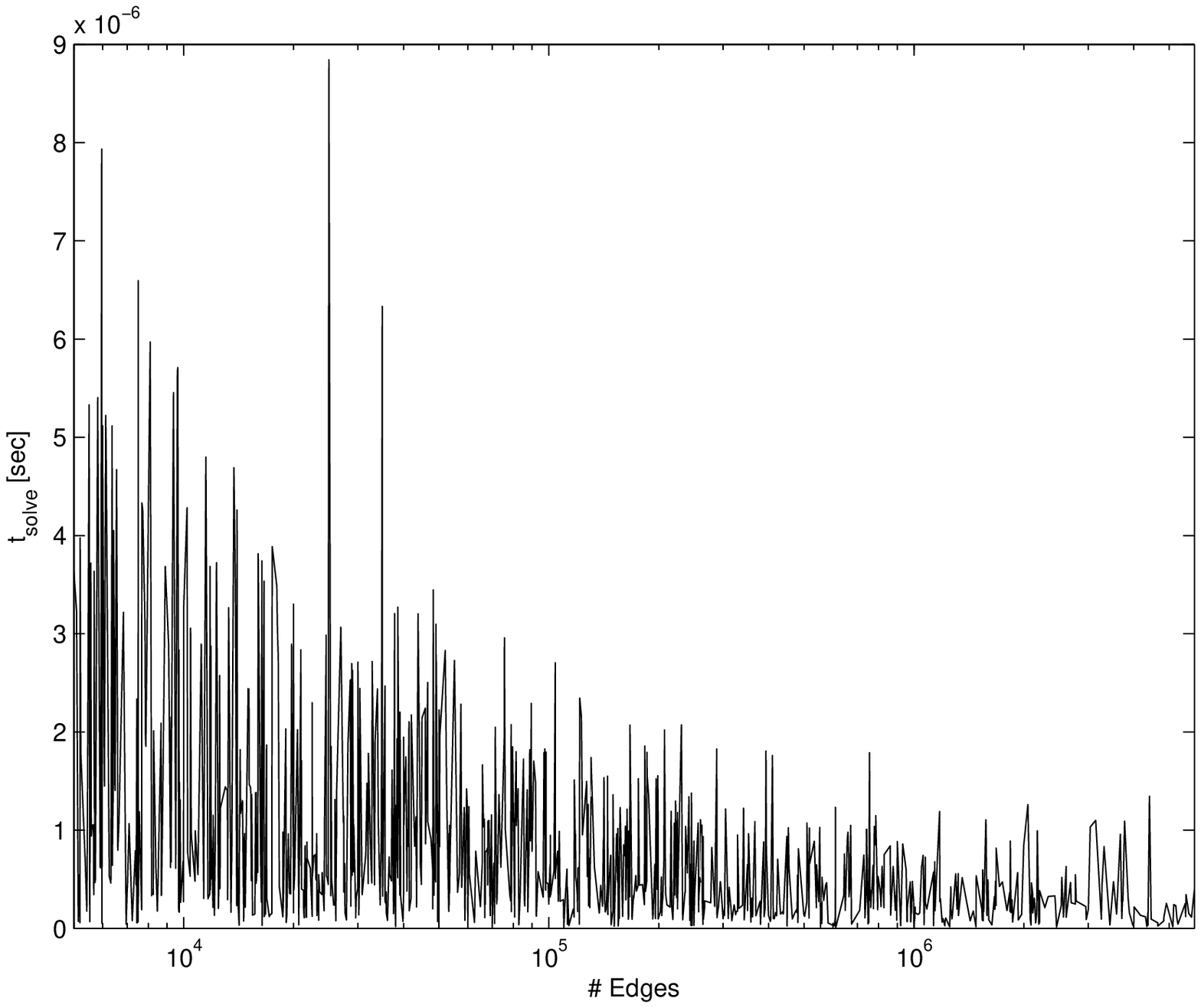}
      \\
      (c) & (d)
    \end{tabular}
    \end{center}
    \caption{LAMG performance vs. number of edges for graphs with $5,000$ edges and larger: (a) Total time per edge $\ttotal$. (b) Adaptive correction solve time gain $A$. The dashed line is the break-even point $A=1$. (c) normalized setup time $\tsetup$. (d) Normalized solve time $\tsolve$.}
    \label{times}
\end{figure}

LAMG scaled linearly with graph size: both $\tsetup$ and $\tsolve$ were approximately constant, and the total time per edge was $2.9 \cdot 10^{-5}$ seconds on average. Stated differently, LAMG performed a single linear solve to $10$ significant figures at $33,000$ edges per second. See Fig.~\ref{times}a,c,d and Table~\ref{times_avg}. 

Adaptive energy correction provided a $15\%-20\%$ speed up and was superior for almost all graphs (cf. Fig.~\ref{times}b). The total time comprised $70\%$ setup and $30\%$ solve for the flat scheme and $75\%/25\%$ for the adaptive scheme. The respective average ACFs were $.18$ and $.048$. These results were better than expected, possibly because of the specific nature of graphs in the collection. We therefore also include details on the typically harder graphs for LAMG in Table~\ref{harder}.

\begin{table}[htbp]
\centering \footnotesize
\begin{tabular}{|l|c|c|c|}
\hline
\multicolumn{2}{|l|}{Measure}     & Median & Mean $\pm$ Std. Dev. \\ \hline \hline
\multicolumn{2}{|l|}{Total time [sec]}           & $2.4 \cdot 10^{-5}$ & $2.9 \cdot 10^{-5} \pm 2.5 \cdot 10^{-5}$ \\ \hline
\multicolumn{2}{|l|}{Setup time [sec]}           & $1.8 \cdot 10^{-5}$ & $2.1 \cdot 10^{-5} \pm 1.5 \cdot 10^{-5}$ \\ \hline
\multirow{2}{*}{Solve time per figure [sec]}  & Flat $\mu=\frac43$ & $6.2 \cdot 10^{-7}$ & $1.3 \cdot 10^{-6} \pm 1.7 \cdot 10^{-6}$ \\ \cline{2-4}
                             & Adaptive  & $5.3 \cdot 10^{-7}$ & $8.7 \cdot 10^{-7} \pm 1.0 \cdot 10^{-6}$ \\ \hline
\multirow{2}{*}{ACF}  & Flat $\mu=\frac43$ & $.180$ & $.239 \pm .209$ \\ \cline{2-4}
                             & Adaptive  & $.048$ & $.069 \pm .068$ \\ \hline
\multicolumn{2}{|l|}{Adaptive gain $A$} & $1.170$ & $1.238 \pm .302$ \\ \hline
\end{tabular}
\caption{Median and mean LAMG performance measures over all graphs with $5,000$ or more edges. ACF statistics only include graphs for which more than one level was constructed.}
\label{times_avg}
\end{table}

\begin{table}[htbp]
\footnotesize
\centering
\begin{tabular}{|p{1.13in}|c|c|c|c|c|c|c|}
\hline
Name&$n\,/\,m$&$M$&$L$&\multicolumn{2}{c|}{ACF}&$\frac{\tsetup}{\ttotal}$&$\ttotal$\\ \cline{5-6}
&&&&\tiny{Flat}&\tiny{Adaptive}&&\\ \hline \hline
Stanford website&$668925 / 3732100$&$5011$&$17$&$.381$&$.114$&$85.9\%$&$1.9 \cdot 10^{-5}$\\ \hline
Calif. road network&$1965206 / 2766607$&$2638$&$19$&$.423$&$.198$&$81.0\%$&$3.5 \cdot 10^{-5}$\\ \hline
GaAsH$_6$ molecule&$61349 / 1660230$&$1$&$9$&$.144$&$.053$&$78.9\%$&$1.6 \cdot 10^{-5}$\\ \hline
Citeseer citations&$268495 / 1156647$&$1$&$13$&$.519$&$.175$&$83.5\%$&$2.7 \cdot 10^{-5}$\\ \hline
Amazon sales&$400727 / 1049624$&$1$&$14$&$.449$&$.172$&$81.2\%$&$3.2 \cdot 10^{-5}$\\ \hline
RF circuit device&$74044 / 207323$&$1$&$14$&$.211$&$.102$&$73.4\%$&$7.8 \cdot 10^{-5}$\\ \hline
Optimization problem&$8364 / 105723$&$1$&$11$&$.173$&$.056$&$72.3\%$&$2.2 \cdot 10^{-5}$\\ \hline
\end{tabular}
\caption{LAMG performance for several graphs representative of the harder cases for LAMG.}
\label{harder}
\end{table}

\subsection{Example of a Hierarchy}
\begin{figure}[htbp]
    \centering
    \begin{center}
    \begin{tabular}{cc}
      \includegraphics[height=1.75in]{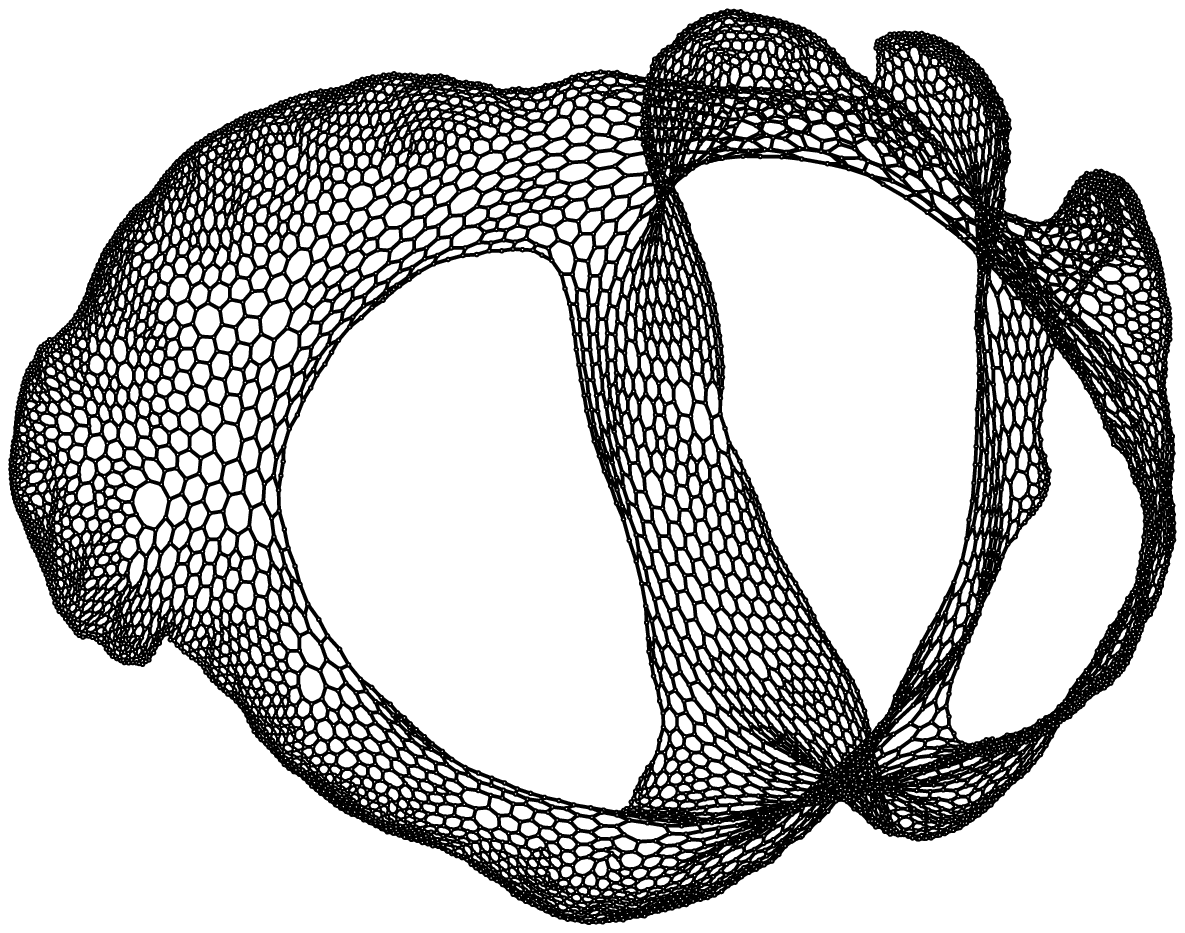}
      &
      \includegraphics[height=1.75in]{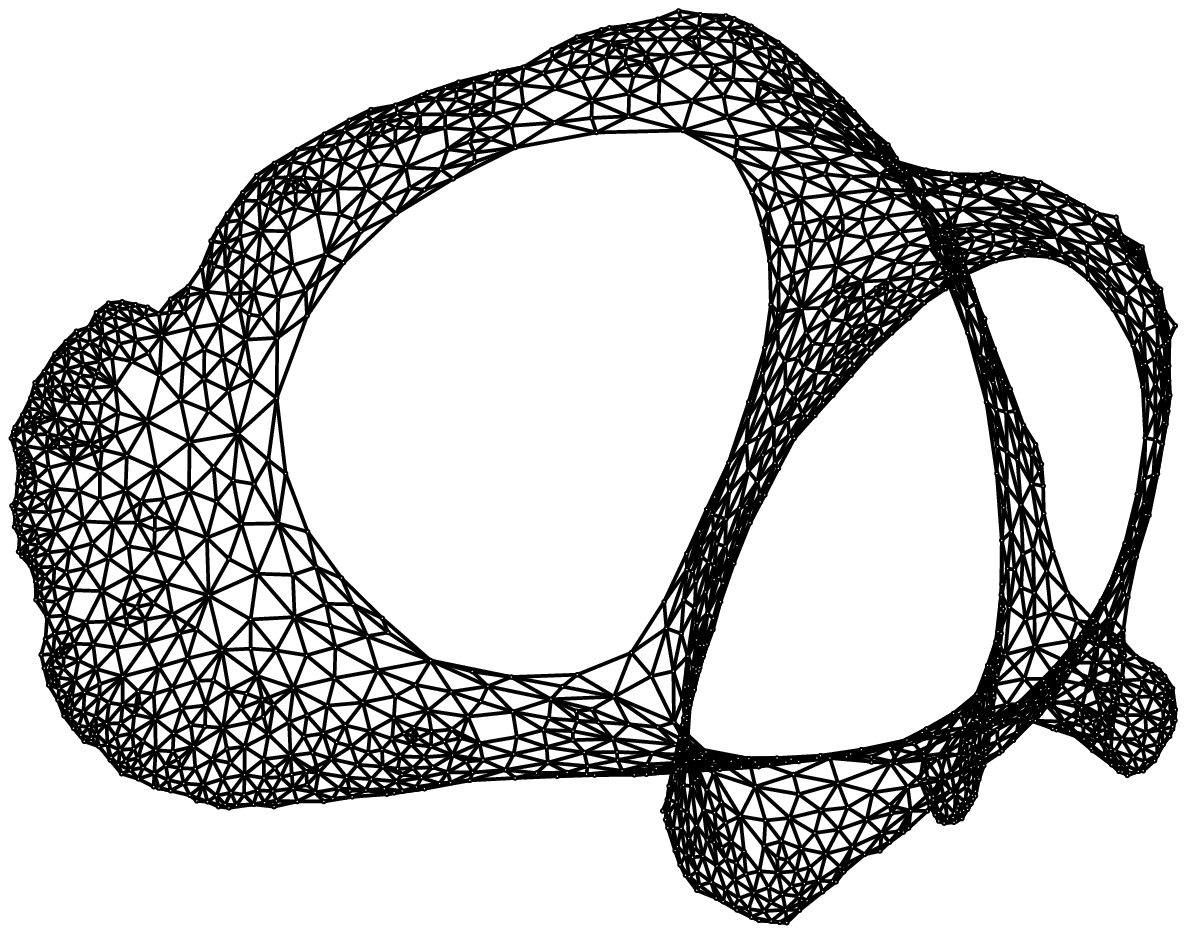}
      \\
      $G^1$ & $G^3$ \\
      \includegraphics[height=1.75in]{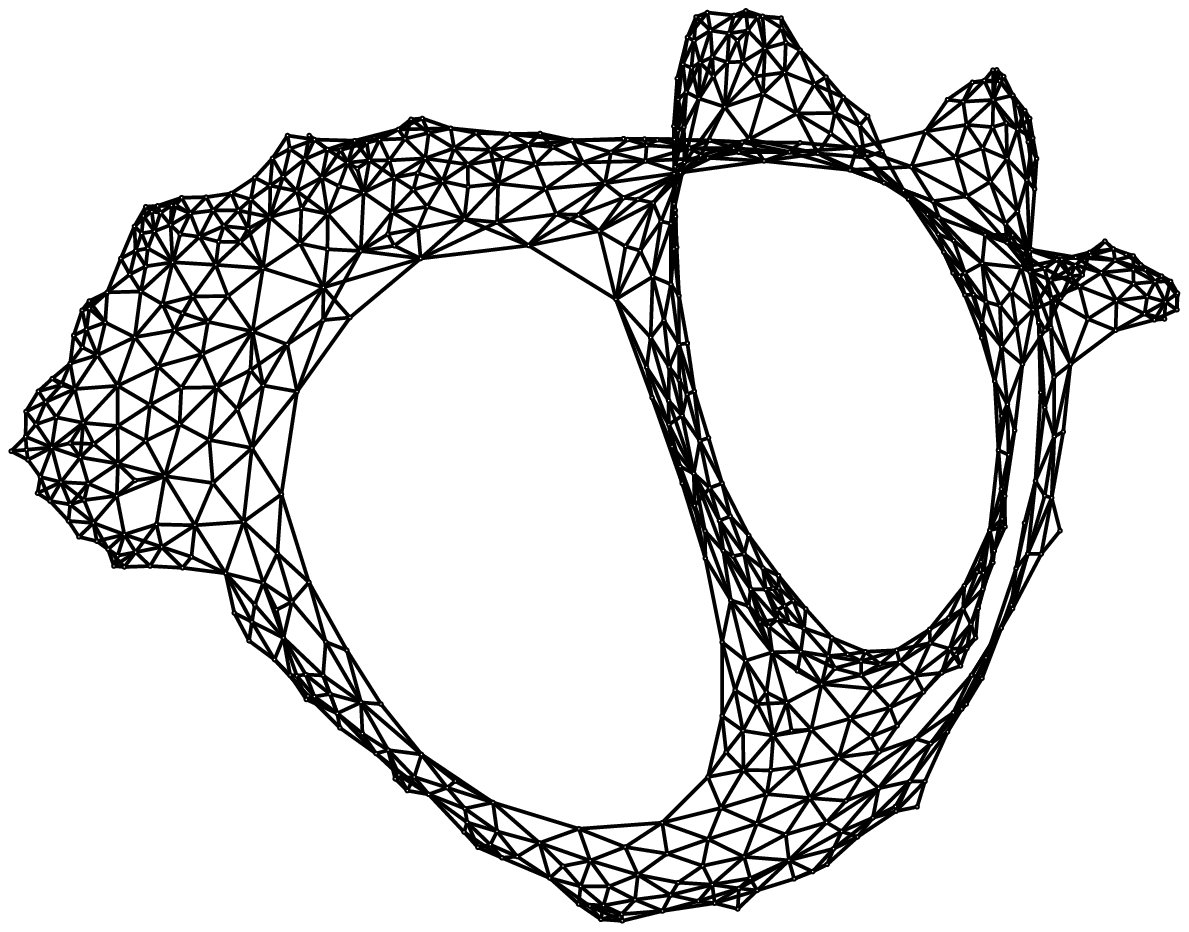}
      &
      \includegraphics[height=1.75in]{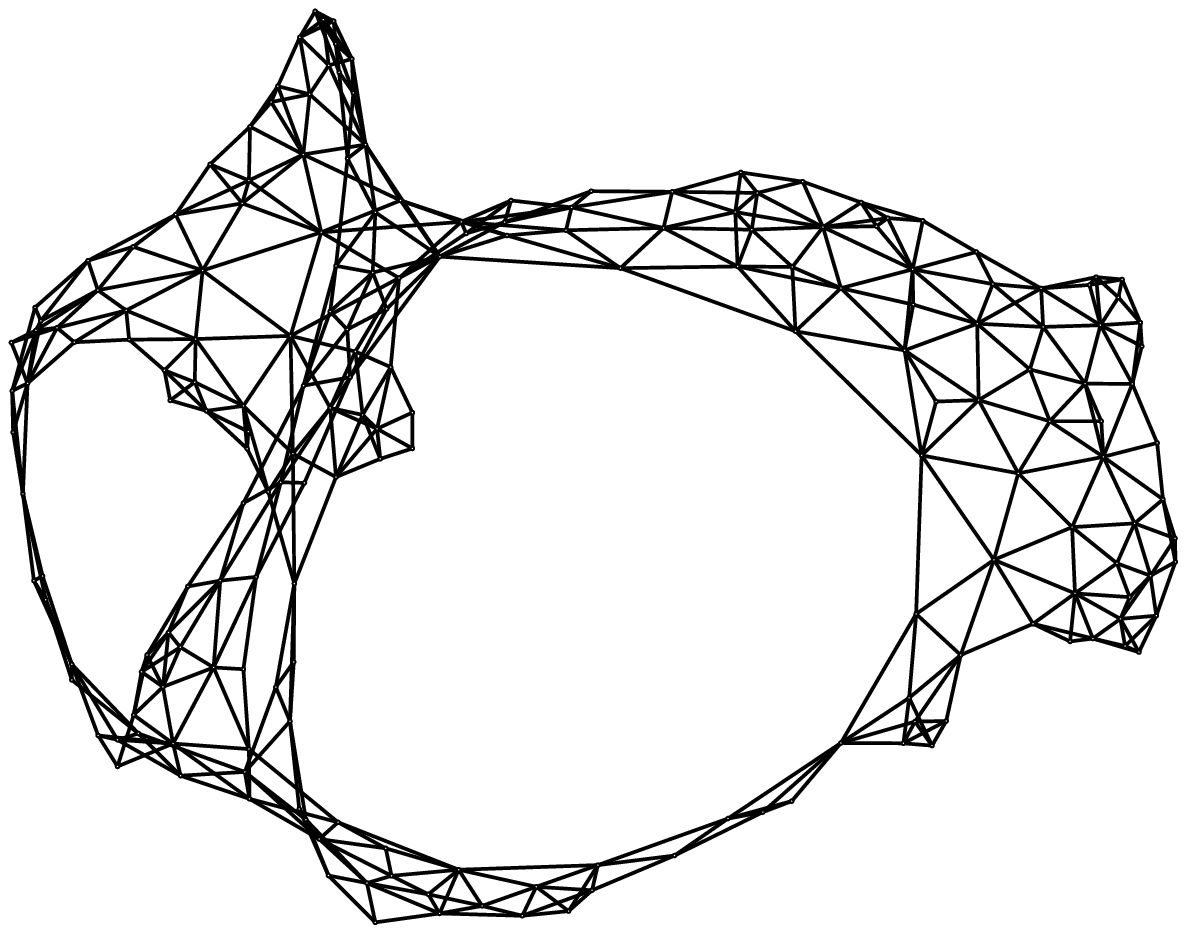}
      \\
      $G^5$ & $G^7$
    \end{tabular}
    \end{center}
    \caption{The four finest non-elimination level graphs for the 2-D airfoil problem.}
    \label{airfoil_graphs}
\end{figure}

For the UF graph {\tt AG-Monien/airfoil1-dual} -- a 2-D air foil finite element singly-connected graph, an eight-level LAMG hierarchy was constructed (see Fig.~\ref{airfoil_graphs} and Table\ref{airfoil_hierarchy}). The total number of edges at all levels (measuring the required storage) was three times the finest graph size. Cycle ACFs were $.358$ and $.157$ for the flat and adaptive energy corrections, respectively.

\begin{table}[htbp]
\footnotesize
\centering
\begin{tabular}{|c|c|c|c|c|c|c|}
\hline
$l$ & Type & $n \,/\, m$ & Mean Degree & $\nu^l_1+\nu^l_2$ & $\gamma^l$ & $K$ \\ \hline \hline
$1$ & \finest &  $8034 / 11813$ & $2.94$ & $   0 $ & $  1.0 $ & $  0 $ \\ \hline
$2$ & \elim   &  $3858 / 11279$ & $5.84$ & $   3 $ & $  1.5 $ & $  8 $ \\ \hline
$3$ & \agg    &  $1561 /  4523$ & $5.80$ & $   0 $ & $  1.0 $ & $  0 $ \\ \hline
$4$ & \elim   &  $1323 /  4059$ & $6.14$ & $   3 $ & $  1.5 $ & $  9 $ \\ \hline
$5$ & \agg    &  $ 615 /  1754$ & $5.70$ & $   0 $ & $  1.0 $ & $  0 $ \\ \hline
$6$ & \elim   &  $ 472 /  1451$ & $6.14$ & $   3 $ & $  1.5 $ & $ 10 $ \\ \hline
$7$ & \agg    &  $ 198 /   530$ & $5.36$ & $   0 $ & $  1.0 $ & $  0 $ \\ \hline
$8$ & \elim   &  $ 114 /   339$ & $5.94$ & $   3 $ & $  2.0 $ & $ 11 $ \\ \hline
\end{tabular}
\caption{LAMG level hierarchy for the 2-D airfoil problem. Graphs were drawn using GraphViz 2.28.0 with the SFDP algorithm \cite{graphviz}.}
\label{airfoil_hierarchy}
\end{table}

\subsection{Grids with Negative Weights}
\label{negative_weights}

Like Bootstrap AMG \cite{cr_etna, alg_distance_anis}, LAMG is not restricted to M-matrices and also applies to some graphs with negative edge weights $w_{uv}$, as long as the Laplacian matrix is (or very close to being) positive semi-definite. This distinguishes LAMG from graph-theoretic works \cite{st06, koutis} and classical AMG. To demonstrate this capability, we tested LAMG on the following 2-D grid Laplacians, whose stencils are depicted in Fig.~\ref{negative_stencils}:
\bi
    \item[(a)] The standard 5-point finite-difference (FD) discretization of $U_{xx} + U_{yy}$ on the unit square with Neumann boundary conditions (BC).
    \item[(b)] The 13-point $\kth{4}$-order FD stencil of $U_{xx} + U_{yy}$.
    \item[(c)] The discretized anisotropic-rotated Laplace operator
    \be
        \left(\cos^2(\alpha) + \ep \sin^2(\alpha)\right) U_{xx} + (1-\ep) \sin(2\alpha) U_{xy}
        + \left(\ep \cos^2(\alpha) + \sin^2(\alpha)\right) U_{yy}\,,
        \label{anis_rot}
    \ee
    with $\alpha=-\pi/4$, $\ep=10^{-4}$, standard 5-point stencil of $U_{xx}$,$U_{yy}$, and an alignment-agnostic cross-term
    $$
        U_{xy} \approx
           \frac{1}{4 h^2}
    \left[
        \begin{tabular}{rrr}
              -1   &    0  &   1   \\
               0   &    0  &   0    \\
               1   &    0  &  -1
        \end{tabular}
    \right]\,,$$
    where $h$ is the grid meshsize. Neumann BC were enforced.

    \item[(d)] The same as (c), but aligning the cross-term with the northeast and southwest neighbors:
    $$
        U_{xy} \approx
           \frac{1}{2 h^2}
    \left[
        \begin{tabular}{rrr}
               0   &   -1  &   1   \\
              -1   &    2  &  -1    \\
               1   &   -1  &   0
        \end{tabular}
        \right]\,.$$
    \item[(e)] The Finite-Element (FE) discretization of the Laplace operator with stretched quadrilateral elements ($h_x/h_y\rightarrow \infty$) with periodic boundary conditions.
    \item[(f)] The 13-point FD discretization of the biharmonic operator $\Delta^2 U$ with the boundary conditions $\partial U/\partial n = \partial^2 U/\partial n^2 = 0$.
\ei
\begin{figure}[htbp]
    \centering\footnotesize
    \begin{center}
    \begin{tabular}{cc}
   $\left[
        \begin{tabular}{rrr}
                   &   -1  &      \\
               -1  &    4  &  -1    \\
                   &   -1  &
        \end{tabular}
    \right]$
    &
    $\left[
        \begin{tabular}{rrrrr}
            &       &   -1  &       & \\
            &       &  -16  &       & \\
        1   & -16   &   60  & -16   & 1 \\
            &       &  -16  &       & \\
            &       &    1  &       &
        \end{tabular}
    \right]$
      \\
      (a) & (b) \\
    $\left[
        \begin{tabular}{rrr}
        -.24998   &   -.50005  &   .24998  \\
        -.50005   &    2.0002  &  -.50005   \\
         .24998   &   -.50005  &  -.24998
        \end{tabular}
    \right]$
    &
    $\left[
        \begin{tabular}{rrr}
                   &   -1  & .49995   \\
              -1   & 3.0001  &  -1     \\
          .49995   &  -1  &
        \end{tabular}
    \right]$
      \\
      (c) & (d) \\
    $\left[
        \begin{tabular}{rrr}
           -1   & -4  &  -1  \\
            2   &  8  &   2  \\
           -1   & -4  &  -1
        \end{tabular}
    \right]$
    &
    $\left[
        \begin{tabular}{rrrrr}
            &       &    1  &       & \\
            &   2   &   -8  &   2   & \\
        1   &  -8   &   20  &  -8   & 1 \\
            &   2   &   -8  &   2   & \\
            &       &    1  &       &
        \end{tabular}
        \right]$
      \\
      (e) & (f) \\
    \end{tabular}
    \end{center}
    \caption{Stencils of grid Laplacians with negative weights. Entries are normalized to a meshsize-independent sum and rounded to five significant figures.}
    \label{negative_stencils}
\end{figure}
Problems (c) and (d) are bad discretizations that do not align with the characteristic direction of (\ref{anis_rot}).
(c),(d) and (e) are considered hard for AMG \cite{alg_distance_anis}.

LAMG exhibited mesh-independent convergence and total time in all cases and scaled linearly with grid size. Performance figures are given in Table~\ref{results_negative}.

\begin{table}[htbp]
\centering\footnotesize
\begin{tabular}{|p{3.5cm}|c|c|c|c|c|c|}
\hline
Problem&$m$&$L$&\multicolumn{2}{c|}{ACF}&$\frac{\tsetup}{\ttotal}$&$\ttotal$\\ \cline{4-5}
&&&\tiny{Flat}&\tiny{Adaptive}&&\\ \hline \hline
(a) 5-point &$523264$&$19$&$0.279$&$0.136$&$71.7\%$&$4.9 \times 10^{-5}$\\ \hline
(b) 13-point $\kth{4}$ order &$1567746$ &$17$&$0.358$&$0.196$&$71.5\%$&$7.4 \times 10^{-5}$\\ \hline
(c) Anis. rot., agnostic &$1045506$&$16$&$0.763$&$0.713$&$45.6\%$&$8.2 \times 10^{-5}$\\ \hline
(d) Anis. rot., misaligned &$784385$&$17$&$0.763$&$0.680$&$49.5\%$&$1.3 \times 10^{-4}$\\ \hline
(e) Stretched FE &$1567746$&$17$&$0.807$&$0.725$&$39.4\%$&$9.3 \times 10^{-5}$\\ \hline
(f) Biharmonic&$1048576$&$14$&$0.789$&$0.731$&$48.3\%$&$7.0 \times 10^{-5}$\\ \hline
\end{tabular}
\caption{LAMG performance for grid graphs on a $512 \times 512$ grid with $n=262144$ nodes.}
\label{results_negative}
\end{table}

While the paper \cite{alg_distance_anis} focused on accurately finding the characteristic directions of (c)--(e) without sparing setup costs and only presented two-level experiments, LAMG is a full multi-level method with a far shorter setup time whose ACF is reasonable in all cases, albeit this ACF can also be significantly reduced using bootstrap tools.

The biharmonic case is intriguing. Even though GS is not the best smoother \cite[\S 20.3.2]{guide} and the interpolation must be second-order \cite[\S 3.3.5]{mg_theory}, the caliber-1-based, black-box LAMG performed reasonably well without any tuning.

These results are certainly only preliminary; further research should be directed toward improving the convergence rates in cases (c)-(f).

\subsection{Lean Geometric Multigrid}
\label{lmg}

Higher performance for the Poisson equation on uniform grid can be obtained by a standard 1:2 coarsening in every dimension at all levels and employing GS Red-Black (RB) relaxation \cite[\S 3.6]{guide}. LAMG then reduces to {\it Lean Geometric Multigrid (LMG)}: standard multigrid \cite[\S 1]{guide} cycle with index $\gamma=1.5$, first-order transfers and energy-corrected coarsening. Since the energy ratio is $2$ for all error modes, we employ a flat correction $\mu=2$ in Eq.~(\ref{galerkin2}).

As a basic experiment, we compared LMG's performance to the ``classical'' cycle with GS-RB, linear interpolation and second-order full weighting in two-level Local Mode Analysis (LMA) and multilevel experiments with seven levels on a $128 \times 128$ grid for the 2-D periodic Poisson problem. The LMG (1,2)-cycle turns out to be a record-breaking Poisson solver in terms of asymptotic efficiency at a $.5$ convergence per unit work, versus the classical V(1,1) at $.67$ (cf. Table~\ref{results_lmg}).

\begin{table}[htbp]
\centering\footnotesize
\begin{tabular}{|c|c||c|c|c|c|}
\hline 
Experiment & Method & \multicolumn{4}{c|}{$(\nu_1,\nu_2)$} \\ \cline{3-6}
&&$(1,0)$ & $(1,1)$ & $(1,2)$ & $(2,2)$ \\ \hline \hline
\multirow{2}{*}{Two-level LMA}  & LMG & $.500\,( 1.5)$ & $.125\,( 2.5)$ & $.034\,( 3.5)   $ & $.024\,( 4.5)$ \\  \cline{2-6}
                                & Classical & $.250\,( 3.3)$ & $.074\,( 4.3)$ & $.051\,( 5.3)       $ & $.039\,( 6.3)$ \\ \hline
\multirow{2}{*}{Multilevel ACF} & LMG ($\gamma=1.5$)       & $.470\,( 2.1)$ & $.119\,( 3.6)$ & $.032\,( 5.0)  $ & $.022\,( 6.4)$ \\ \cline{2-6}
& Classical V-cycle	     & $.261\,( 4.4)$ & $.101\,( 5.7)$ & $.060\,( 7.1) 	  $ & $.045\,( 8.4)$ \\ \hline
\end{tabular}
\caption{LMA predictions and actual asymptotic convergence of LMG vs. classical multigrid for the 2-D periodic Poisson problem. The ACF and cycle work in relaxation units (in parentheses) are listed for each case.}
\label{results_lmg}
\end{table}

This caliber-1 ``super-efficiency'' for the {\it periodic} Poisson problem can be explained by interpreting the grids of the hierarchy as cell-centered discretizations at meshsize $h,2h,4h,\dots$, for which the LMG restriction operator is effectively second-order when $\mu=2$. Although the interpolation is still first-order, LMG still satisfies the transfer operator rules for optimal multigrid efficiency \cite[\S 4.3]{guide}.

Other boundary conditions impede the choice of $\mu$. While supplementary local relaxations near boundaries theoretically ensure attaining the two-level rates, it would be more beneficial to study the performance of adaptive energy correction in LMG.

\section{Extensions}
\label{extensions}
Enhancements and adaptations of the LAMG approach to related computational problems are outlined.

\subsection{High Performance Implementation}

The main bottleneck of the current \matlab implementation is the aggregation, specifically, locating the neighbors $\cA_u$ of node $u$ and the corresponding graph weights, needed for the various term computations in Algorithm~\ref{bestSeed2}. This is because the aggregation stage sequentially scans and updates nodes while \matlab excels at vector operations. An informal survey of language operation benchmarks seems to indicate that a C or C++ implementation could provide a 2--4 speed-up. The $mex$ interface can be used to program these internal loops in C/C++ while continuing to use \matlab for the entire program \cite{matlab_mex}.

AMG parallelization to multiple processors is nontrivial. Lessons from other parallel AMG works \cite{fischer, blatt} will likely directly apply to LAMG. Regarding LAMG-specific operations, TVs and affinities can be computed in parallel; parallelizing the aggregation stage will again pose the greatest challenge. Aggregation decisions may need to be modified to be symmetric to avoid conflicts near processor boundaries. An important advantage is LAMG's tendency to reduce the size of the coarse stencils, which should help reduce the overlap between the sub-graphs assigned to neighboring processors, as well as communication at the coarsest levels.

We plan to compare LAMG's performance within existing AMG infrastructures such as Hypre and Trilinos-PETSc. Like classical AMG solvers, LAMG can be incorporated as either a stand-alone solver or as a CG preconditioner.

\subsection{Coarsening Improvements}
\label{improvements}
The final word has by no means been said on the details of the LAMG algorithm of \S \ref{algorithm}.
\bi
    \item Currently, $u$ can only be aggregated with a directly-neighboring seed $s$. In some problems, one should also search within $u$'s $\bA^2$-neighbors to construct a good aggregation. For instance, in the anisotropic-rotated problem Fig.~\ref{negative_stencils}d, $u$ should be aggregated along the characteristic direction, i.e., with its southeast (or northwest) neighbor, neither of which is contained in $\cA$.
    \item During the aggregation stage (Algorithm~\ref{aggregationStage}), a clever ordering of the undecided nodes can winnow out ``holes'', i.e., lone nodes that cannot be associated just because all their neighbors have been aggregated, and no aggregate enlargement can be warranted by the energy ratio control.
    \item If no small energy ratio can be found, or if subsequent cycle convergence is slow, isolated bottleneck nodes may be de-aggregated.
    \item Bigotry of the ``lean'' approach should not be practiced, either: one can occasionally up the interpolation caliber at these troublesome nodes, provided that this does not substantially increase the total coarse edges.
    \item Since the coarsening algorithm utilizes both affinity and energy ratio conditions, it may be possible to reduce the number of TVs $K$ and/or the number of relaxation passes $\nu$ (or trade one for the other to increase efficiency at a given complexity).
\ei

\subsection{Local Energy Corrections}
\label{other_corrections}
Instead of a flat $\mu=\frac43$ factor in (\ref{galerkin2}), one can apply different $\mu$'s to different aggregates. We experimented with different energy correction schemes, some based on fitting the coarse nodal energies of TVs to their fine counterparts (Eq.~(\ref{nodal_coarse})). While this can dramatically curtail energy inflation, care must be taken to avert over-fitting that ultimately results in the coarse-level correction operator's instability.

Analogously, one can define a local adaptive $\mu$ in the MINRES procedure of \S \ref{adaptive} at each level $l$. For example, $\mu$ may be constant on the nodes of each level $l+1$ aggregate. $\mu$ should be smoothed by (say) a GS relaxation sweep on $\bA^{l+1} \bm{\mu} = \bzero$, as it multiplies a smooth correction vector and cannot be allowed to radically oscillate.

It is unclear whether the extra work would be justified in either case; it may be useful for problems re-solved for many RHS vectors or when a larger setup overhead is tolerable.

\subsection{The Eigenvalue Problem}
\label{eigenvalue}

AMG can be nicely combined with the Exact Interpolation Scheme (EIS) \cite{eis} to find the smallest nonzero eigenvalue and associated Fiedler vector. After the LAMG setup phase, multilevel EIS cycles can be applied to $\bA \bx = \lambda \bx$ using the sparsity pattern of $\bP$ and $\bA$ at all levels. After the current approximation $\tbx$ is relaxed at level $l$, the interpolation weights are re-set to $p^{l+1}_{uv} = \tilde{x}_u/\tilde{x}_v$ (and possibly modified back to $1$ near $\tilde{x}$'s zeros), followed by recomputing $\bA^{l+1} = (\bP^l)^{T} \bA^l \bP^l$ as well as coarsening the mass matrix $\bB^{l+1}=(\bP^l)^{T} \bB^l \bP^l$; $\bB^1 := \bI$. At the coarsest level, the lowest eigenpair of the pencil $(\bA^L,\bB^L)$ is calculated. The same $\gamma, \nu_1, \nu_2$ can be used as in the linear solver.

The reasoning behind this algorithm is that a piecewise-constant interpolation fits {\it all} near-null-space eigenvectors of $\bA$, not just the constant vector. Hence, large affinities indicate large correlations among the nodal values of the Fiedler vector, so the coarsening {\it pattern} is adequate for its computation. EIS is particularly attractive thanks its super-linear convergence.

Alternatively, one can incorporate the LAMG linear solver into a Rayleigh Quotient iteration, which converges locally cubically to the smallest eigenvalue \cite[\S 8.2]{gvl}.

We plan on developing the LAMG eigensolver (including its generalization to finding several lowest eigenpairs) in the near future, because the Laplacian eigenproblem is even more ubiquitous in applications than the linear system.

\subsection{Other Linear Systems}
\label{other}

The caliber-1 algorithm can be applied to non-zero row sum matrices, except that the interpolation weights are no longer $1$. The affinity definition (\ref{cuv}) remains intact, and the corresponding $\bP$ entry is set to
\be
    p_{uv} := \argmin{q} \left\|X_u-q X_v\right\|_{uv}^2
\ee
(cf. (\ref{affinity_set}). Normally, relaxed TVs yield an accurate enough $p_{uv}$; in problems with almost-zero modes, e.g., the QCD gauge Laplacian, TVs may need to be improved by a bootstrap cycle \cite{bamg}.

Further research should be conducted for negative-weight graphs such as the high-order finite element and anisotropic grid graphs of \S \ref{negative_weights}. The reported convergence factors can be reduced by producing bootstrapped TVs via applying multilevel cycles to $\bA \bx=\bzero$. The cycle is far more powerful than plain relaxation in damping smooth characteristic components, which should lead to more meaningful algebraic distances and to the correct anisotropic coarsening in a second setup round (much larger spacing in the characteristic direction and no coarsening in the cross-characteristic direction). The bootstrap procedure should be similarly useful for many other graphs.

\section{Conclusion}
Over the last decade, Laplacian matrices have attracted increasing attention because they underlie a plethora of graph computational applications ranging from genetic data clustering to social networks to fluid dynamics. To the best of our knowledge, the presented algorithm, Lean Algebraic Multigrid (LAMG), is the first graph Laplacian linear solver demonstrated to scale linearly with graph size. LAMG will hopefully pave the way to significant speed-ups in those applications.

\bibliographystyle{siam}       
\bibliography{lamg}

\end{document}

%% file: styles.tex
\usepackage{graphicx}
\usepackage{algorithm}
\usepackage{algorithmic}
\usepackage{amsmath}
\usepackage{amsfonts}
\usepackage{amssymb}
\usepackage{bm}
\usepackage{url}
\usepackage{multirow}

\newcommand{\ep}{\varepsilon}

\newcommand{\Etot}{E_{\mathrm{tot}}}
\newcommand{\Imax}{I_{\mathrm{max}}}
\newcommand{\numgraphs}{1666\;}
\newcommand{\maxedges}{six million\;}

\newcommand{\by}{\bm{\mathrm{y}}}
\newcommand{\bx}{\bm{\mathrm{x}}}
\newcommand{\bb}{\bm{\mathrm{b}}}
\newcommand{\bu}{\bm{\mathrm{u}}}
\newcommand{\bq}{\bm{\mathrm{q}}}
\newcommand{\bee}{\bm{\mathrm{e}}}
\newcommand{\br}{\bm{\mathrm{r}}}
\newcommand{\balpha}{\bm{\mathrm{\alpha}}}
\newcommand{\tbx}{\tilde{\bx}}

\newcommand{\tbee}{\tilde{\bee}}
\newcommand{\bzero}{\bm{\mathrm{0}}}
\newcommand{\bxi}{\bm{\mathrm{\xi}}}

\newcommand{\bA}{\bm{\mathrm{A}}}
\newcommand{\bB}{\bm{\mathrm{B}}}
\newcommand{\bC}{\bm{\mathrm{C}}}

\newcommand{\bI}{\bm{\mathrm{I}}}
\newcommand{\bL}{\bm{\mathrm{L}}}
\newcommand{\bN}{\bm{\mathrm{N}}}
\newcommand{\bU}{\bm{\mathrm{U}}}
\newcommand{\bP}{\bm{\mathrm{P}}}
\newcommand{\bQ}{\bm{\mathrm{Q}}}
\newcommand{\bR}{\bm{\mathrm{R}}}
\newcommand{\bT}{\bm{\mathrm{T}}}
\newcommand{\bX}{\bm{\mathrm{X}}}
\newcommand{\bPi}{\bm{\mathrm{\Pi}}}

\newcommand{\cN}{\mathcal{N}}
\newcommand{\cE}{\mathcal{E}}
\newcommand{\cA}{\mathcal{A}}
\newcommand{\cT}{\mathcal{T}}
\newcommand{\cZ}{\mathcal{Z}}
\newcommand{\cF}{\mathcal{F}}
\newcommand{\cC}{\mathcal{C}}
\newcommand{\cU}{\mathcal{U}}
\newcommand{\cV}{\mathcal{V}}

\newcommand{\Real}{\mathbb{R}}

\newcommand{\matlab}{\textsc{Matlab\;}}
\newcommand{\finest}{\textsc{Finest\;}}
\newcommand{\elimnospace}{\textsc{Elimination}}
\newcommand{\elim}{\textsc{Elimination\;}}
\newcommand{\agg}{\textsc{Aggregation\;}}
\newcommand{\notvisited}{\textsc{NotVisited}}
\newcommand{\fnode}{\textsc{FNode}}
\newcommand{\notelim}{\textsc{NotEliminated}}
\newcommand{\undecided}{\textsc{Undecided}}
\newcommand{\notfound}{\textsc{NotFound}}
\newcommand{\seed}{\textsc{Seed}}
\newcommand{\visited}{visited}
\newcommand{\status}{status}
\newcommand{\aggsize}{aggregateSize}
\newcommand{\nodes}{\mathrm{nodes}}
\newcommand{\tsetup}{t_{\mathrm{setup}}}
\newcommand{\tsolve}{t_{\mathrm{solve}}}
\newcommand{\ttotal}{t_{\mathrm{total}}}
\newcommand{\sz}{\mathrm{size}}

\newcommand{\kth}[1]{#1^\mathrm{th}}
\newcommand{\lp}{\left(}
\newcommand{\rp}{\right)}

\newcommand{\loget}{\log(1/\ep)}

\newcommand{\mymin}[1]{\underset{#1}{\operatorname{min}}\;}
\newcommand{\mymax}[1]{\underset{#1}{\operatorname{max}}\;}
\newcommand{\argmin}[1]{\underset{#1}{\operatorname{argmin}}\;}
\newcommand{\argmax}[1]{\underset{#1}{\operatorname{argmax}}\;}
\newcommand{\amax}{\alpha_{\mathrm{max}}}

\newcommand{\be}{\begin{equation}}
\newcommand{\ee}{\end{equation}}
\newcommand{\bes}{\begin{equation*}}
\newcommand{\ees}{\end{equation*}}
\newcommand{\bea}{\begin{eqnarray}}
\newcommand{\eea}{\end{eqnarray}}
\newcommand{\beas}{\begin{eqnarray*}}
\newcommand{\eeas}{\end{eqnarray*}}
\newcommand{\bi}{\begin{itemize}}
\newcommand{\ei}{\end{itemize}}
\newcommand{\ben}{\begin{enumerate}}
\newcommand{\een}{\end{enumerate}}

%% file: lamg.bbl
\begin{thebibliography}{10}

\bibitem{add}
{\sc P.~R. Amestoy, T.~A. Davis, and I.~S. Duff}, {\em An approximate minimum
  degree ordering algorithm}, SIAM J. Mat. Anal. Appl., 17 (1996),
  pp.~886--905.

\bibitem{templates}
{\sc R.~Barrett, M.~Berry, T.~F. Chan, J.~Demmel, J.~Donato, J.~Dongarra,
  V.~Eijkhout, R.~Pozo, C.~Romine, and H.~Van der Vorst}, {\em Templates for
  the Solution of Linear Systems: Building Blocks for Iterative Methods, 2nd
  Edition}, SIAM, Philadelphia, PA, 1994.

\bibitem{blatt}
{\sc M.~Blatt}, {\em A parallel algebraic multigrid method for elliptic
  problems with highly discontinuous coefficients}, PhD thesis, Universit{\"a}t
  Heidelberg, 2010.

\bibitem{boman}
{\sc E.~G. Boman, B.~Hendrickson, and S.~Vavasis}, {\em Solving elliptic finite
  element systems in near-linear time with support preconditioners}, SIAM J.
  Num. Anal., 46 (2008), pp.~3264--3284.

\bibitem{poly_book}
{\sc P.~Borwein and T~Erd{\'e}lyi}, {\em Polynomials and Polynomial
  Inequalities}, Springer-Verlag, New York, 1995.

\bibitem{braess}
{\sc D.~Braess}, {\em Towards algebraic multigrid for elliptic problems of
  second order}, Computing, 55 (1995), pp.~379--393.

\bibitem{mg_theory}
{\sc A.~Brandt}, {\em Rigorous quantitative analysis of multigrid, {I.}
  {C}onstant coefficients two-level cycle with $l_2$-norm}, Siam J. Num Anal.,
  31 (1994).

\bibitem{cr_etna}
\leavevmode\vrule height 2pt depth -1.6pt width 23pt, {\em General highly
  accurate algebraic coarsening}, J. Electron. Trans. Num. Anal., 10 (2000),
  pp.~1--20.
\newblock Multilevel methods (Copper Mountain, CO, 1999).

\bibitem{yes}
\leavevmode\vrule height 2pt depth -1.6pt width 23pt, {\em Multiscale
  scientific computation: Review 2001}, in Multiscale and Multiresolution
  Methods, T.~Barth, T.~Chan, and R.~Haimes, eds., Springer-Verlag, 2002,
  pp.~3--96.

\bibitem{alg_distance_anis}
{\sc A.~{Brandt}, J.~{Brannick}, K.~{Kahl}, and I.~{Livshits}}, {\em {An
  algebraic distances measure of AMG strength of connection}}, ArXiv e-prints,
  (2011).

\bibitem{bamg}
{\sc A.~Brandt, J.~Brannick, K.~Kahl, and I.~Livshits}, {\em Bootstrap amg},
  SIAM J. Sci. Comp., 33 (2011), pp.~612--632.

\bibitem{guide}
{\sc A.~Brandt and O.~E. Livne}, {\em Multigrid Techniques: 1984 Guide with
  Applications to Fluid Dynamics}, Classics in Applied Mathematics, SIAM,
  revised~ed., 2011.

\bibitem{geodetic}
{\sc A.~Brandt, S.~McCormick, and J.~W. Ruge}, {\em Algebraic multigrid ({AMG})
  for automatic multigrid solution with application to geodetic computations},
  tech. report, Colorado State University, Fort Collins, Colorado, 1983.

\bibitem{cr_james}
{\sc J.~J. Brannick and R.~D. Falgout}, {\em Compatible relaxation and
  coarsening in algebraic multigrid}, SIAM J. Sci. Comp., 32 (2010),
  pp.~1393--1416.

\bibitem{asa}
{\sc M.~Brezina, R.~Falgout, S.~MacLachlan, T.~Manteuffel, S.~McCormick, and
  J.~Ruge}, {\em Adaptive smoothed aggregation ({ASA})}, SIAM J. Sci. Comp., 25
  (2004), p.~2004.

\bibitem{sa}
{\sc M.~Brezina, P.~Vanek, and P.~S. Vassilevski}, {\em An improved convergence
  analysis of smoothed aggregation algebraic multigrid}, Num. Lin. Alg. Appl.,
  (2011).

\bibitem{hongkong_class}
{\sc H.~Chang and D.~Yeung}, {\em Graph laplacian kernels for object
  classification from a single example}, in In CVPR (2, 2006, pp.~2011--2016.

\bibitem{chung}
{\sc F.~R.~K. Chung}, {\em {Spectral Graph Theory (CBMS Regional Conference
  Series in Mathematics, No. 92)}}, American Mathematical Society, Feb. 1997.

\bibitem{DS08}
{\sc S.~I. Daitch and D.~A. Spielman}, {\em Faster approximate lossy
  generalized flow via interior point algorithms}, CoRR, abs/0803.0988 (2008).

\bibitem{uf_collection}
{\sc T.~A. Davis}, {\em University of florida sparse matrix collection}, NA
  Digest, 92 (1994).

\bibitem{dimacs}
{\sc {DIMACS Center at Rutgers Univesity}}, {\em Dimacs implementation
  challenges}.
\newblock Available online at \url{http://dimacs.rutgers.edu/Challenges/},
  2011.

\bibitem{ding}
{\sc C.~H.~Q. Ding, X.~He, H.~Zha, M.~Gu, and H.~D. Simon}, {\em A min-max cut
  algorithm for graph partitioning and data clustering}, in Proceedings of ICDM
  2001, 2001, pp.~107--114.

\bibitem{regression}
{\sc N.R. Draper and H.~Smith}, {\em Applied Regression Analysis},
  Wiley-Interscience, 1998.

\bibitem{graphviz}
{\sc J.~Ellson, E.~R. Gansner, E.~Koutsofios, S.~C. North, and G.~Woodhull},
  {\em Graphviz - open source graph drawing tools}, Graph Drawing,  (2001),
  pp.~483--484.

\bibitem{hypre}
{\sc R.~Falgout, A.~Cleary, J.~Jones, E.~Chow, V.~Henson, C.~Baldwin, P.~Brown,
  P.~Vassilevski, and U.~Meier Yang}, {\em Hypre reference manual 2.7.0b}.
\newblock
  \url{https://computation.llnl.gov/casc/hypre/download/hypre-2.7.0b_ref_manual.pdf},
  2011.

\bibitem{fischer}
{\sc P.~Fischer, J.~Lottes, D.~Pointer, and A.~Siegel}, {\em Petascale
  algorithms for reactor hydrodynamics}, J. Physics: Conference Series, 125
  (2008), p.~012076.

\bibitem{FG07}
{\sc A.~Frangioni and C.~Gentile}, {\em Prim-based support-graph
  preconditioners for min-cost flow problems}, Comput. Optim. Appl., 36 (2007),
  pp.~271--287.

\bibitem{matlab_sparse}
{\sc J.~R. Gilbert, C.~Moler, and R.~Schreiber}, {\em Sparse matrices in
  matlab: Design and implementation}, SIAM J. Mat. Anal. Appl., 13 (1992),
  pp.~333--356.

\bibitem{gvl}
{\sc G.H. Golub and C.F.V. Loan}, {\em {Matrix Computations}}, Johns Hopkins
  studies in the mathematical sciences, Johns Hopkins University Press,
  third~ed., 1996.

\bibitem{class04}
{\sc L.~Gorelick, M.~Galun, E.~Sharon, R.~Basri, and A.~Brandt}, {\em Shape
  representation and classification using the poisson equation}, in In In Proc.
  of CVPR'04, 2004, pp.~61--67.

\bibitem{trilinos}
{\sc M.~A. Heroux, R.~A. Bartlett, V.~E. Howle, R.~J. Hoekstra, J.~J. Hu, T.~G.
  Kolda, R.~B. Lehoucq, K.~R. Long, R.~P. Pawlowski, E.~T. Phipps, A.~G.
  Salinger, H.~K. Thornquist, R.~S. Tuminaro, J.~M. Willenbring, A.~Williams,
  and K.~S. Stanley}, {\em An overview of the trilinos project}, ACM Trans.
  Math. Softw., 31 (2005), pp.~397--423.

\bibitem{koutis}
{\sc I.~Koutis, G.~L. Miller, and R.~Peng}, {\em Approaching optimality for
  solving sdd linear systems}, in Proceedings of the 2010 IEEE 51st Annual
  Symposium on Foundations of Computer Science, FOCS '10, Washington, DC, USA,
  2010, IEEE Computer Society, pp.~235--244.

\bibitem{clustering06}
{\sc D.~Kushnir, M.~Galun, and A.~Brandt}, {\em Fast multiscale clustering and
  manifold identification}, Pattern Recognition, 39 (2006), pp.~1876--1891.

\bibitem{eis}
\leavevmode\vrule height 2pt depth -1.6pt width 23pt, {\em Efficient multilevel
  eigensolvers with applications to data analysis tasks}, IEEE Transactions on
  Pattern Analysis and Machine Intelligence, 32 (2010), pp.~1377--1391.

\bibitem{ancestry}
{\sc Ann~B. Lee, Diana Luca, Lambertus Klei, Bernie Devlin, and Kathryn
  Roeder}, {\em Discovering genetic ancestry using spectral graph theory},
  Genet. Epidem., 34 (2010), pp.~51--59.

\bibitem{cr_oren}
{\sc O.~E. Livne}, {\em Coarsening by compatible relaxation}, Num. Lin. Alg.
  Appl., 11 (2004), pp.~205--227.

\bibitem{wave_accuracy}
{\sc I.~Livshits and A.~Brandt}, {\em Accuracy properties of the wave-ray
  multigrid algorithm for helmholtz equations}, SIAM J. on Sci. Comp., 28
  (2006), pp.~1228--1251.

\bibitem{matlab_mex}
{\sc {MathWorks}}, {\em Creating {C/C++} language mex-files}.
\newblock Available online at
  \url{http://www.mathworks.com/help/techdoc/matlab_external/f7667.html}, 2011.

\bibitem{olson}
{\sc L.~N. Olson, J.~B. Schroder, and R.~S. Tuminaro}, {\em A general
  interpolation strategy for algebraic multigrid using energy minimization},
  SIAM J. Sci. Comp., 33 (2011), pp.~966--991.

\bibitem{alg_distance}
{\sc D.~Ron, I.~Safro, and A.~Brandt}, {\em Relaxation-based coarsening and
  multiscale graph organization}, Multiscale Model. Sim., 9 (2011),
  pp.~407--423.

\bibitem{rs86}
{\sc J.~Ruge and K.~St{\"u}ben}, {\em Algebraic multigrid ({AMG})}, in
  Multigrid Methods, Frontiers in Applied Mathematics, S.~F. McCormick, ed.,
  SIAM, 1987, pp.~73--130.

\bibitem{mloga}
{\sc I.~Safro}, {\em Minimum logarithmic arrangement ({MinLogA}) results
  archive}.
\newblock Available online at \url{http://www.mcs.anl.gov/~safro/mloga.html},
  2011.

\bibitem{sedgewick}
{\sc R.~Sedgewick}, {\em Algorithms in {C++}, Part 5: Graph Algorithms},
  Addison-Wesley, 2002.

\bibitem{radu}
{\sc A.~Sharma, R.~P. Horaud, D.~Knossow, and E.~von Lavante}, {\em Mesh
  segmentation using laplacian eigenvectors and gaussian mixtures}, in
  Proceedings of AAAI Fall Symposium on Manifold Learning and its Applications,
  Fall Symposium Series Technical Reports, Arlington, VA, November 2009, AAAI
  Press.

\bibitem{nature}
{\sc E.~Sharon, M.~Galun, D.~Sharon, R.~Basri, and A.~Brandt}, {\em Hierarchy
  and adaptivity in segmenting visual scenes}, Nature, 442 (2006),
  pp.~810--813.

\bibitem{icm10}
{\sc D.~A. Spielman}, {\em Algorithms, graph theory, and linear equations in
  laplacian matrices}, in Proceedings of the International Congress of
  Mathematicians 2010 (ICM 2010), World Scientific, 2010, pp.~2698--2722.

\bibitem{st06}
{\sc D.~A. Spielman and S.~Teng}, {\em Nearly-linear time algorithms for
  preconditioning and solving symmetric, diagonally dominant linear systems},
  CoRR, abs/cs/0607105 (2006).

\bibitem{tarjan}
{\sc R.~E. Tarjan}, {\em Depth first search and linear graph algorithms}, SIAM
  J. Computing, 1 (1972), pp.~146--160.

\bibitem{md}
{\sc W.~F. Tinney and J.~W. Walker}, {\em {Direct solutions of sparse network
  equations by optimally ordered triangular factorization}}, Proc. IEEE, 55
  (1967), pp.~1801--1809.

\bibitem{trot}
{\sc U.~Trottenberg, C.~W. Oosterlee, and A.~Sch{\"u}ller}, {\em Multigrid},
  Academic Press, London, 2000.

\bibitem{walshaw}
{\sc C.~Walshaw}, {\em Multilevel refinement for combinatorial optimisation
  problems}, Annals Oper. Res., 131 (2004), pp.~325--372.

\bibitem{schur_book}
{\sc F.~Zhang}, {\em The Schur Complement And Its Applications (numerical
  Methods And Algorithms)}, Springer-Verlag, 2005.

\bibitem{zhu_learning}
{\sc X.~Zhu, Z.~Ghahramani, and J.~D. Lafferty}, {\em Semi-supervised learning
  using {G}aussian fields and harmonic functions.}, in ICML'03, 2003,
  pp.~912--919.

\end{thebibliography}
